\documentclass{tac}
\copyrightyear{2007}

\title{Cartesian Bicategories II}
\author{A. Carboni, G.M. Kelly, R.F.C Walters, and R.J. Wood}
\thanks{The authors gratefully acknowledge financial support from 
the Australian ARC, and the Canadian NSERC.
Diagrams typeset using M. Barr's diagram package, diagxy.tex.}
\address{\\Dipartimento di Scienze delle Cultura\\
Politiche e dell'Informazione\\
Universit\`a dell Insubria, Italy\\[3pt] 
and\\[3pt]
School of Mathematics\\
University of Sydney\\
NSW 2006, Australia\\[3pt]
and\\[3pt]
\\Dipartimento di Scienze delle Cultura\\
Politiche e dell'Informazione\\
Universit\`a dell Insubria, Italy\\[3pt]
and\\[3pt]
Department of Mathematics and Statistics\\
Dalhousie University\\
Halifax, NS, B3H 3J5 Canada}

\eaddress{aurelio.carboni@uninsubria.it, maxk@usyd.edu.au,\\ 
robert.walters@uninsubria.it, rjwood@mathstat.dal.ca}

\newtheorem{axiom}{Axiom}

\let\thm\theorem
\let\prp\proposition
\let\cor\corollary
\let\lem\lemma
\let\dfn\definition
\let\eth\endtheorem
\let\prf\proof
\let\frp\endproof
\let\exl\example
\let\rmk\remark

\input diagxy
\xyoption{curve}

\newcommand{\s}{\scalefactor{0.5}}

\newcommand{\f}{{\kern -.25em}:{\kern -.25em}}

\newcommand{\ra}{{\s\to}}
\newcommand{\la}{\,{\s\toleft}\,}

\newcommand{\rass}{\s{|{\kern -.38em}\to}}

\newcommand{\iso}{\cong}
\newcommand{\laj}{\dashv}

\newcommand{\x}{\times}
\newcommand{\ox}{\otimes}

\newcommand{\op}{{^\mathrm{op}}}

\newcommand{\inv}{^{-1}}

\newcommand{\CAT}{\mathbf{CAT}}

\newcommand{\one}{\mathbf{1}}

\newcommand{\cart}{\mathrm{Cart}}

\newcommand{\map}{\mathrm{Map}}

\newcommand{\prof}{\mathrm{Prof}}

\newcommand{\spn}{\mathrm{Span}}

\newcommand{\E}{{\cal E}}

\newcommand{\bA}{\mathbf{A}}
\newcommand{\bB}{\mathbf{B}}

\newcommand{\bG}{\mathbf{G}}
\newcommand{\bM}{\mathbf{M}}

\begin{document}

\maketitle

\begin{abstract}
The notion of {\em cartesian bicategory\/}, introduced in
\cite{caw} for locally ordered bicategories, is extended
to general bicategories. It is shown that a cartesian bicategory
is a symmetric monoidal bicategory.
\end{abstract}

\section{Introduction}\label{intro}
\subsection{}\label{int}
We recall that in \cite{caw} a locally ordered bicategory
$\bB$ was said to be {\em cartesian\/} if the subbicategory of left
adjoints, $\map\bB$, had finite products; each hom-category
$\bB(B,C)$ had finite products; and a certain derived tensor
product on $\bB$, extending the product structure of $\map\bB$,
was functorial. It was shown that cartesian structure
provides an elegant base for sets of axioms characterizing
bicategories of
\begin{enumerate}
\item[i)] relations in a regular category
\item[ii)] ordered objects and order ideals in an exact category
\item[iii)] additive relations in an abelian category
\item [iv)] relations in a Grothendieck topos.
\end{enumerate}
Notable was an axiom, {\em groupoidalness\/}, that
captures the {\em discrete\/} objects in a cartesian locally ordered
bicategory and gives rise to a very satisfactory approach to duals.

It was predicted in \cite{caw} that the notion of cartesian
bicategory would be developable without the restriction of local
orderedness, so as to capture
\begin{enumerate}
\item[v)] spans in a category with finite limits
\item[vi)] profunctors in an elementary topos.
\end{enumerate}
It is this development of the unrestricted notion of cartesian
bicategory that is our present concern. In a sequel \cite{ckww} 
we shall give such further characterization theorems but this
paper is concerned with the basic development. 

A cartesian bicategory is a bicategory with various {\em properties\/}. 
In contrast to many bicategorical studies, {\em no constraint data
subject to coherence conditions are assumed\/}. Constraint data is
constructed from the existence clauses of the universal properties 
provided by the definition of cartesian bicategory and it is tempting 
to believe that all coherence conditions for such constraints will 
automatically follow from the uniqueness clauses of such universal
properties. After all, a mere category with finite products has
$-\x-$ as a candidate for a monoidal structure $-\ox-$. In this
case it is easy to construct arrows 
$a_{X,Y,Z}\f(X\x Y)\x Z\ra X\x(Y\x Z)$ and prove that they are
invertible and satisfy Mac~Lane's pentagon condition, for example.

A cartesian bicategory need not have all finite products (even in the
sense of bilimits). However, it admits a canonical $-\ox-$ connective
that one immediately suspects to be the main ingredient of a
symmetric monoidal structure. One of the axioms for a cartesian
bicategory requires that the locally full subbicategory determined by
the left adjoint arrows admits all finite products in the sense of
bilimits. We show that another bicategory associated with a
cartesian bicategory also has finite products in this sense.
Assuming that a bicategory with finite products underlies a
symmetric monoidal bicategory we are then able to show that a
cartesian bicategory also underlies a symmetric monoidal bicategory.

When we started this project we assumed it was well known that
a bicategory with finite products in the sense of bilimits {\em is\/}
a symmetric monoidal bicategory. All to whom we spoke about the matter  
agreed with us but we were unable to find a published account.
The universal property for products in the sense of bilimits is 
given in terms of equivalences (rather than isomorphisms) of
hom-categories. At first this seems to be a rather weak beginning
when one considers that a monoidal bicategory is a one-object
tricategory whose definition appears very complicated to the 
uninitiated. It seemed that rather a lot had to be constructed,
and even more to be proved, starting from very little. To make 
matters worse we need symmetry in our sequel
paper \cite{ckww} and a symmetric monoidal bicategory is a
one-object, one-arrow, one-2-cell, one-3-cell weak 6-category.
Of course one has the more informative definitions of a
symmetric monoidal bicategory provided collectively by \cite{das} 
and \cite{McC} but we ultimately decided that our assumptions about
finite products required proofs.

Accordingly, we begin in Section 2 with the study of a bicategory $\bA$ 
with finite products, concluding it with Theorem \ref{finprodimpmonl}
stating that $(\bA,\x,1,\cdots)$ is a symmetric monoidal bicategeory. In
Section 3 we find it convenient to define a bicategory $\bB$ to
be {\em precartesian} if the locally full subbicategory determined
by the left adjoints, $\bM:=\map\bB$, has finite products and all
the hom categories have finite products (as mere categories). We
construct from $\bB$ and $\bM$ a further bicategory $\bG$, whose
objects are general arrows of $\bB$ and whose arrows between them
are given by squares containing a 2-cell in which the arrow
components are left adjoints. We show that $\bG$ too has finite
products, for a precartesian $\bB$, and use $\bG$ to define
arrows 
$$\ox\f\bB\x\bB\ra\bB\la\one\f I$$
which are proved to be lax functors. Finally, in Section 4 we
define $\bB$ to be cartesian if it is precartesian and enjoys
the further property that the constructed constraints for $\ox$
and $I$ are invertible, thus making these pseudofunctors. We
prove Theorem \ref{smbicat}, that $(\bB,\ox,I,\cdots)$ is
a symmetric monoidal bicategeory, and further basic facts that will
be used in \cite{ckww}.

G.M. ``Max'' Kelly died during the preparation of this paper. He
was actively working on it on the day of his passing. The other 
authors express their gratitude for his work here and for so much more
that he had shared with us as a friend and a colleague over many years.
We regret too that he was unable to provide a proof reading of our
final draft.

\section{Bicategories with Finite Products}\label{finprod}
\subsection{}\label{notations}
We consider a bicategory $\bA$.
We write $X,Y,Z,...$ for objects,
$R,S,T,...$ for arrows, and $\alpha,\beta,\gamma,...$ for 
2-cells of $\bA$. 
We usually omit parentheses in three-fold composites but then
$RST$ is to be understood as $(RS)T$. In this vein we define
$\bA(T,R)\f\bA(X,A)\ra\bA(Y,B)$ as the 
composite 
$\bA(X,A)\to^{\bA(1,R)}\bA(X,B)\to^{\bA(T,1)}\bA(Y,B)$,
where $R\f A\ra B$ and $T\f Y\ra X$ in $\bA$, 
and use this choice in defining the hom pseudofunctor
\begin{equation}\label{hom}
\bA(-,-)\f\bA\op\x\bA\ra\CAT
\end{equation}
where $\CAT$ denotes the 2-category of categories. Note that we use
{\em pseudofunctor} for what is also called a {\em homomorphism
of bicategories} and {\em lax functor} for what is also
called a {\em morphism of bicategories}.

We find it convenient to assume all our bicategories to be
{\em normal}, in the sense that the constraints $1R\cong R$
and $R1\cong R$ are identities. It is easy to replace any
bicategory by a biequivalent normal one: we have only to
provide a new identity.

A functor $F\f X\ra Y$, or more generally an arrow $F\f X\ra Y$ in
a bicategory, is said to be an {\em equivalence} if there is an
arrow $U\f Y\ra X$ and
invertible 2-cells $\epsilon\f FU\ra 1$ and $\alpha\f 1\ra UF$. It is
a bicategorical formality that in such a situation one can
find an invertible $\eta\f 1\ra UF$ giving an adjunction 
$\eta,\epsilon\f F\laj U$. 
We may as well suppose each equivalence to come with a 
specified right adjoint with invertible unit and invertible counit.

An arrow $R$ in a bicategory is called a {\em map\/} if it has a 
right adjoint. 
We shall usually denote maps by lowercase Roman letters; and 
if $f$ is a map, we write $\eta_f,\epsilon_f\f f\laj f^*$ for a chosen 
adjunction that makes it so.
We write $\map\bA$ for the locally-full 
subbicategory of $\bA$ determined by the maps.

\subsection{}\label{prt} 
We begin by supposing the bicategory $\bA$ to have binary and nullary 
{\em products} (in the {\em bilimit} sense, 
which is the only one appropriate to bicategories). 
For the binary case, this means simply that to each 
pair $(X,Y)$ of objects is assigned an object $X \x Y$ and arrows 
$p_{X,Y}\f X\x Y\ra X$ and $r_{X,Y}\f X\x Y\ra Y$, called projections,
such that the functor $\bA(A, X\x Y) \ra \bA(A,X) \x\bA(A,Y)$
they determine is an equivalence for each $A$. (We use lowercase
for projections because in the case of interest they are hypothesized
to be maps.)
In elementary terms, $\bA(A, X\x Y) \ra \bA(A,X) \x\bA(A,Y)$
is essentially surjective on objects and fully faithful. Thus for
$R\f A \ra X$ and $S\f A\ra Y$ we have
$\langle R,S\rangle\f A\ra X\x Y$
and invertible 2-cells
$\mu_{R,S}\f  p_{X,Y}\langle R,S\rangle\ra R$ and
$\nu_{R,S}\f  r_{X,Y}\langle R,S\rangle\ra S$. And for
$T,U\f A\ra X\x Y$ together with 2-cells 
$\alpha\f p_{X,Y}T\ra p_{X,Y}U$ and $\beta\f r_{X,Y}T\ra r_{X,Y}U$ 
there exists a unique $\gamma\f T\ra U$ having 
$p_{X,Y}\gamma=\alpha$ and $r_{X,Y}\gamma=\beta$. In particular, for
2-cells $\phi\f R\ra R'\f A\ra X$ and $\psi\f S\ra S'\f A\ra Y$,
there is a unique 2-cell 
$\langle\phi,\psi\rangle\f\langle R,S\rangle\ra\langle R',S'\rangle$
making $\mu_{R,S}$ and $\nu_{R,S}$ natural in $R$ and $S$.

Next, since {\em each} pair of objects has such a product, we can, 
as is well known, make $\x$ into a pseudofunctor 
$\x\f\bA\x\bA\ra\bA$ 
in such a way that the $p_{X,Y}$ and $r_{X,Y}$ become the components 
of pseudonatural transformations $p$ and $r$. In more detail, for
$R\f X\ra A$ and $S\f Y\ra B$ we have $R\x S\f X\x Y\ra A\x B$ and 
invertible 2-cells $p_{R,S}\f p_{A,B}(R\x S)\ra R p_{X,Y}$ and 
$r_{R,S}\f r_{A,B}(R\x S)\ra S r_{X,Y}$, making $p$ and $r$ into 
pseudonatural transformations $p\f \x \ra P$ and $r\f \x\ra R$, 
where $P,R\f\bA\x\bA\ra\bA$ are the projection pseudofunctors. 
Of course we take $R\x S$  to be 
$\langle R p_{X,Y}, S r_{X,Y}\rangle$ in the language of the last 
paragraph, while $p_{R,S}$ is a suitable $\mu$
and $r_{R,S}$ is a suitable $\nu$; but we don't have to remember 
their origins.  

As a matter of fact we have a problem here with notation, because 
later we use $p_{R,S}$ in a different sense in connection with
the product projection $R\x S\ra R$ in another bicategory $\bG$ whose
objects are the arrows of a bicategory $\bB$. 
To avoid confusion we henceforth denote the 2-cells of the last 
paragraph, namely the
pseudonaturality isomorphisms of $p$ and $r$, by $p'_{R,S}$
and $r'_{R,S}$ instead of $p_{R,S}$ and $r_{R,S}$. More generally,
if $t$ is a pseudonatural transformation with components $t_X$, we
will write $t'_R$ for the pseudonaturality isomorphism corresponding
to an arrow $R\f X\ra Y$.

Of course $\x\f\bA\x\bA\ra\bA$ is the right pseudoadjoint of the 
diagonal pseudofunctor $\Delta\f\bA\ra\bA\x\bA$, with 
counit constituted by the pseudonatural $p$ and $r$. 
For the unit $d$ we can take the components $d_X$ to be 
$\langle1,1\rangle\f X\ra X\x X$, with the pseudonaturality 
isomorphism $d'_R\f d_YR\ra(R\x R)d_X$ coming from the evident 
isomorphisms of its domain and its codomain with $\langle R,R\rangle$.

There remains the nullary product, namely the {\em terminal object} 
of $\bA$. Here we are given an object $1$ of $\bA$ for which the 
functor $\bA(A,1)\ra \one$ is an equivalence for each A, so that 
$1\f \one\ra\bA$ is the right pseudoadjoint of $!\f \bA\ra \one$. The 
unit is given by choosing an arrow $t_X\f X\ra 1$  for each $X$; 
then for $R\f X\ra A$ there is a unique 2-cell $t'_R\f t_AR\ra t_X$, 
and it is invertible.

\subsection{}\label{n>-0} 
The bicategory $\bA$ is said to have {\em finite products} if, for
every finite set $I$ and every family $X=(X_i)_{i\in I}$ of 
objects of $\bA$, there is an object $P$ and a family of arrows 
$p=(p_i\f P\ra X_i)_{i\in I}$ with domain $P$ so that, for each $A$, 
the functor
\begin{equation}\label{one}
\bA(A,P)\ra\prod_{i\in I}\bA(A,X_i)              
\end{equation}
induced by the $p_i$ is an equivalence. Such a family 
$p=(p_i\f P\ra X_i)_{i\in I}$ is called a {\em product cone over}~$X$.  

\prp\label{0&2=>fin} 
A bicategory with binary products and a terminal object has
finite products.
\eth
\prf
The cases of $I$ having cardinality $2$ or $0$ are covered explicitly 
by the hypotheses. They are a matter of identifying $\bA^I$ with
$\bA\x\bA$ or $\one$. For $I$ having cardinality $1$ and an 
$I$-indexed family $X$, any equivalence $P\ra X$, in 
particular the identity $1_X\f X\ra X$, is a product cone.
Suppose that we have a product cone for each $I$-indexed family $X$
with cardinality of $I$ less than $n+1$ and consider a family
$Z$ indexed by a set $J$ of cardinality $n+1$. We can see $J$
as a sum of sets $I+\one$ so that the $J$-family $Z$ is an 
$I$-family $X$ together with a single object $Y$. Let
$(p_i\f P\ra X_i)_{i\in I}$ be a product cone for $X$ and
consider the $J$-family of arrows $q=(q_j\f P\x Y\ra Z_j)_{j\in J}$ 
given by $q_j=p_i.p_{P,Y}\f P\x Y\ra X_i$, for $j=i\in I$ and
$q_j=r_{P,Y}\f P\x Y\ra Y)$ for $j=*\in\one$.
Since a composite of equivalences is an equivalence, this family is
a product cone over $Z$.
\frp

Because a finite set $I$ with cardinality $n$ admits a bijection
to the set $\{1,2,\cdots,n\}$, the bicategory $\bA$ with binary
and nullary products, equivalently all finite products, admits,
for each natural number $n$, a pseudofunctor
$$\Pi_n\f\bA\x\cdots\x\bA=\bA^n\ra\bA$$
with domain the $n$-fold product bicategory, right pseudoadjoint
to the diagonal pseudofunctor. We can write, for example, 
$\Pi_3(X,Y,Z)=X\x Y\x Z$ {\em without parentheses}.

\subsection{}\label{goal} 
We want to exhibit $\bA$ as underlying a monoidal bicategory with 
$\x=\Pi_2$ as its tensor product. To this end, for a family 
$X=(X_i)_{i\in I}$ we write $\bA(X) = \bA((X_i)_{i\in I})$ for the 
bicategory whose objects are the product cones over $X$, an arrow 
from $(b_i\f B\ra  X_i)$ to $(c_i\f C\ra X_i)$ being an
arrow $R\f B\ra C$ in $\bA$ along with isomorphisms 
$\mu_i\f c_iR\ra b_i$, and a 2-cell $(R,\mu_i)\ra(S,\nu_i)$ being
a 2-cell $\alpha\f R\ra S$ in $\bA$ having $\nu_i.(c_i\alpha)=\mu_i$. 
There is an evident forgetful pseudofunctor from $\bA(X)$ to $\bA$.

But since we know that finite products exist, we can 
abbreviate in the above by writing an object $(b_i\f B\ra X_i)$ of 
$\bA(X)$ as $b\f B\ra X$, treating $X$ here as a name of {\em a} 
product  $\prod_{i\in I}X_i$, so that $b\f B\ra X$ is an equivalence. 
Now an arrow $(R,\mu)\f (B,b)\ra(C,c)$ of $\bA(X)$ is an arrow 
$R\f B\ra C$ of $\bA$ along with a single isomorphism $\mu\f cR\ra b$,
and a 2-cell $(R,\mu)\ra(S,\nu)$ is an $\alpha\f R\ra S$ satisfying 
$\nu.(c\alpha) = \mu$.

Recall that a bicategory is biequivalent to the bicategory $\one$
precisely when
\begin{enumerate}
\item[i)] the set of  objects is not empty,
\item[ii)] for any objects B and C there is an arrow $R\f B\ra C$, and
\item[iii)] for any two arrows $R,S\f B\ra C$ there is a unique 
2-cell $R\ra S$.
\end{enumerate}

\prp\label{3.2} 
Each bicategory $\bA(X)$ is biequivalent to $\one$.
\eth
\prf
First, the existence of finite products in $\bA$ ensures that each 
$\bA(X)$ is not empty. Next, for any objects $(B,b)$ and $(C,c)$ of
$\bA(X)$, there is an arrow $(R,\mu)$ from $(B,b)$ to $(C,c)$ because,
$(C,c)$ being a product cone, (\ref{one}) (with $P$ replaced by $C$) 
is essentially surjective. Finally, if $(S,\nu)$ is a second arrow 
from $(B,b)$ to $(C,c)$, then to give a 2-cell 
$\alpha\f (R,\mu)\ra (S,\nu)$ is to give an $\alpha\f R\ra S$ 
having $c_i\alpha = (\nu_i)^{-1}(\mu_i)$; and there is a unique such 
2-cell because, $(C,c)$ being a product cone, (\ref{one}) 
(with $P$ replaced by $C$) is fully faithful.
\frp

Note that, as a consequence, every arrow in $\bA(X)$ is an 
equivalence; so that its underlying arrow in $\bA$ is also an 
equivalence. Again, every 2-cell in $\bA(X)$ is an isomorphism, 
so that its underlying 2-cell in $\bA$ is also an isomorphism.

\subsection{}\label{3.3} 
In the definition of a monoidal bicategory the tensor product 
(here $\x=\Pi_2$) is to be a pseudofunctor. We have remarked in 
\ref{prt} that $\x=\Pi_2$ is a pseudofunctor $\bA\x \bA\ra\bA$, 
and that $p$, $r$, and $t$ are pseudonatural. 
Then there is a pseudofunctor 
$\Pi_2\x 1\f\bA\x\bA\x\bA\ra\bA\x\bA$, and hence a pseudofunctor 
$\Pi_2(\Pi_2\x1)\f\bA\x\bA\x\bA\ra\bA$ sending $(X,Y,Z)$ to 
$(X\x Y)\x Z$. 
Similarly there is a pseudofunctor 
$\Pi_2(1\x \Pi_2)\f\bA\x\bA\x\bA\ra\bA$ 
sending $(X,Y,Z)$ to $X\x(Y\x Z)$. Each of the pseudofunctors 
$\Pi_2(\Pi_2\x 1)$ and $\Pi_2(1\x \Pi_2)$ is an object of 
the bicategory $[\bA^3,\bA]$
of pseudofunctors, pseudonatural transformations, and modifications.

Starting with the pseudonatural $p$ and $r$ we construct the 
projections 
\begin{equation}
p_{X,Y}.p_{(X\x Y),Z},\quad r_{X,Y}.p_{(X\x Y),Z},\quad\mbox{and}\quad 
r_{(X\x Y),Z}     \label{three}
\end{equation}
from $(X\x Y)\x Z$ to $X$, $Y$, and $Z$ respectively, and the 
projections
\begin{equation}
p_{X,(Y\x Z)},\quad p_{Y,Z}.r_{X,(Y\x Z)},\quad\mbox{and}\quad
r_{Y,Z}.r_ {X,(Y\x Z)}   \label{four}     
\end{equation}
from $X\x(Y\x Z)$ to $X$, $Y$, and $Z$ respectively. Like $p$ and $r$,
the projections (\ref{three}) and (\ref{four}) are pseudonatural 
when we treat $(X\x Y)\x Z$, $X\x(Y\x Z)$, $X$, $Y$, and $Z$ as
objects of $[\bA^3,\bA]$, regarding them as alternative names for 
$\Pi_2(\Pi_2\x 1)$, for $\Pi_2(1\x \Pi_2)$, and for the first, 
second, and third projections of $\bA^3$ onto $\bA$. 
Accordingly, we can see the 
projections (\ref{three}) and (\ref{four}) as arrows of $[\bA^3,\bA]$.

The projections (\ref{three}) and (\ref{four}) constitute product 
cones, exhibiting $(X\x Y)\x Z$ and $X\x (Y\x Z)$ as products
in $\bA$ of $X$, $Y$, and $Z$. However products in $[\bA^3,\bA]$ are 
formed pointwise from those in $\bA$. Accordingly, when we see the 
projections (\ref{three}) and (\ref{four}) as arrows of $[\bA^3,\bA]$, 
they form product cones there, exhibiting $(X\x Y)\x Z$ and 
$X\x(Y\x Z)$ as products of $X$, $Y$, and $Z$ not only in $\bA$ but 
also in $[\bA^3,\bA]$.

\subsection{}\label{3.4} 
We now apply Proposition \ref{3.2} to the bicategory 
$[\bA^3,\bA](X,Y,Z)$. We have $\Pi_3\f\bA^3\ra\bA$. The product cones 
(\ref{three}) and (\ref{four}) correspond to equivalences
$h\f \Pi_2(\Pi_2\x1)\ra \Pi_3$ and $k\f \Pi_2(1\x \Pi_2)\ra \Pi_3$, 
with components $h_{X,Y,Z}\f (X\x Y)\x Z\ra X\x Y\x Z$  and 
$k_{X,Y,Z}\f X\x(Y\x Z)\ra X\x Y\x Z$.

It follows from Proposition \ref{3.2} that there is a pseudonatural 
transformation
\begin{equation}
a\f \Pi_2(\Pi_2\x1)\ra \Pi_2(1\x \Pi_2)  \label{five} 
\end{equation}
and an invertible modification
\begin{equation}
\mu\f ka\ra h             \label{six}
\end{equation}
Moreover $a$ is an equivalence, and $(a,\mu)$ is unique to within 
a unique isomorphism. The components of the pseudonatural $a$ are 
equivalences
$$a_{X,Y,Z}\f (X\x Y)\x Z\ra X\x (Y\x Z)$$
in $\bA$, and the components of $\mu$ are isomorphisms
$$\mu_{X,Y,Z}\f k_{X,Y,Z}.a_{X,Y,Z}\ra h_{X,Y,Z}$$

In a similar way we produce the pseudonatural equivalences 
$l\f\Pi_2(\Pi_0\x 1)\ra\Pi_1$ and $r\f\Pi_1\ra\Pi_2(1\x\Pi_0)$.

\subsection{}\label{3.5} 
In $[\bA^4,\bA]$ we have (writing $XY$ for $X\x Y$ and so on) the 
composite
\begin{equation}
((XY)Z)W\ra (X(YZ))W\ra X((YZ)W)\ra X(Y(ZW))    \label{seven} 
\end{equation}
of $a_{X,Y,Z}W$, $a_{X,YZ,W}$, and $Xa_{Y,Z,W}$; on the face of it 
these are arrows in $\bA$, but we can also see them (as in the second 
paragraph of \ref{3.3}) as naming pseudonatural equivalences which 
are arrows of $[\bA^4,\bA]$. We also have the
composite
\begin{equation}
((XY)Z)W\ra(XY)(ZW)\ra X(Y(ZW))                    \label{eight}
\end{equation}
of  $a_{XY,Z,W}$ and $a_{X,Y,ZW}$. Let us write (\ref{seven}) and
(\ref{eight}) as $m,n\f U\ra V$,
where $U=\Pi_2(\Pi_2(\Pi_2\x1)\x1)$ and $V=\Pi_2(1\x \Pi_2(1\x \Pi_2))$;
and denote by $u\f U\ra \Pi_4$ and $v\f V\ra \Pi_4$ the equivalences 
formed from projections like (\ref{three}) and (\ref{four}).

Each of these arrows in (\ref{seven}) and (\ref{eight}) comes with an 
invertible modification, namely $\mu_{X,Y,Z}W$, $\mu_{X,YZ,W}$, 
$X\mu_{Y,Z,W}$, $\mu_{XY,Z,W}$, and $\mu_{X,Y,ZW}$ respectively, 
where $\mu$ is the modification in (\ref{six}).
These modifications fit together and compose to give invertible 
modifications $\alpha\f vm\ra u$ and $\beta\f vn\ra u$. 
Now $(m,\alpha)$ and $(n, \beta)$ are arrows from $U$ to $V$ in 
$[\bA^4,\bA](X,Y,Z,W)$, so that 
there is by Proposition \ref{3.2} a unique 2-cell 
$\pi\f(m,\alpha)\ra(n,\beta)$, which is invertible. That is to say, 
$\pi\f m\ra n$ is an invertible modification satisfying  
$\beta.(v\pi)=\alpha$.

This is the modification $\pi$ of (TD7) in \cite{gps}. There
are three similar modifications called $\mu$, $\lambda$, and $\rho$ in 
\cite{gps} involving respectively $(a,l,r)$, $(a,l,l)$ and $(a,r,r)$, 
and these appear similarly as 2-cells in $[\bA^3, \bA](X,Y,Z)$.

\subsection{}\label{3.6} 
On page 10 of \cite{gps} there is a diagram (TA1) to be satisfied in a
monoidal bicategory. It demands the equality in $[\bA^5,\bA]$ of 
two modifications
$$X(Ya).Xa.a.(Xa)V.aV.(aU)V\two a.a.a$$ 
They are in fact equal by Proposition \ref{3.2}, as they are 2-cells 
in $[\bA^5,\bA](X,Y,Z,U,V)$.

The remaining two axioms follow in a similar way.

\subsection{}\label{3.7} 
We next show that the monoidal bicategory we have constructed is 
symmetric; or more correctly that we can endow it with a symmetry. 
According to \cite{das}, a {\em symmetry} for a monoidal bicategory 
consists of a braiding and a syllepsis, with the syllepsis satisfying 
a certain symmetry condition: see Definition 18 on page 131 
of \cite{das}. Although \cite{das} consider only the special case of 
a Gray monoid, where the associativity is an identity, their 
definition above of a symmetry is surely meant to apply generally. 
The meanings of braiding and syllepsis (but not of symmetry) are 
given in \cite{McC}, as follows (see pages~133 to 145).

The basic datum for a braiding is a pseudonatural equivalence
\begin{equation}
s\f\ox\ra\ox S                    \label{nine}
\end{equation}
where $S\f\bA^2\ra\bA^2$ sends $(X,Y)$ to $(Y,X)$ and so on, 
and where $\ox\f\bA^2\ra\bA$ is the tensor product, which for us is 
$\x$, here more conveniently called $\Pi$. There are two further data, 
consisting of invertible 
modifications sitting in hexagonal diagrams, and the data are to 
satisfy four axioms which are equations between invertible 
modifications.

A syllepsis for a braided monoidal bicategory is an invertible 
modification $\sigma$ from the identity of $\ox$ to 
$sS.s\f\ox\ra\ox SS=\ox$, with components
\begin{equation}
\sigma_{X,Y}\f1_{\ox}\ra s_{Y,X}s_{X,Y}       \label{ten}               
\end{equation}
which is to satisfy two axioms consisting of equations between 
invertible modifications. Note that $1_{\ox}S=1_{\ox S}$ and we have
also $\sigma S\f 1_{\ox S}\ra s.sS$, with components
$$\sigma_{Y,X}\f1_{\ox S}\ra s_{X,Y}s_{Y,X}$$
The braiding and the syllepsis constitute 
a symmetry if (see \cite{das} p.131]) the syllepsis satisfies the 
further condition $s\sigma=(\sigma S)s$, which in terms of
components is
\begin{equation}
s_{X,Y}\sigma_{X,Y}=\sigma_{Y,X}s_{X,Y}\label{eleven} 
\end{equation}

\subsection{}\label{3.8} 
For the monoidal bicategory arising as above from a bicategory $\bA$
with finite products, we construct a symmetry as follows. 
In $[\bA^2,\bA]$ we have the product cone
$$P\toleft^p\Pi\to^rR$$ 
of \ref{prt}, where $P,R\f\bA^2\ra\bA$ are the projection 
pseudofunctors, but also
$$P=RS\toleft^{rS}\Pi S\to^{pS}PS=R$$
is a product cone. By the essential surjectivity aspect of the 
universal property of the latter, we have in the top two triangles 
on the left below an arrow $s$ (necessarily a pseudonatural
equivalence), and invertible 2-cells $\mu$ and $\nu$ as shown.
\begin{equation}\label{twelve}
\bfig
\dtriangle(0,0)|lmm|/->`->`<-/[\Pi`P`\Pi S;p`s`rS]
\btriangle(500,0)|mrm|[\Pi`\Pi S`R;s`r`pS]
\qtriangle(0,-500)|mlm|/<-`<-`->/[P`\Pi S`\Pi;rS`p`sS]
\ptriangle(500,-500)|mmr|/->`->`<-/[\Pi S`R`\Pi;pS`sS`r]
\morphism(300,250)|b|<125,-125>[`;\mu]
\morphism(700,250)|b|<-125,-125>[`;\nu]
\morphism(300,-50)|b|<125,-125>[`;\nu S]
\morphism(700,-50)|b|<-125,-125>[`;\mu S]
\place(1250,0)[=]
\dtriangle(1500,0)|lmm|/->`->`/[\Pi`P`\Pi S;p`s`]
\btriangle(2000,0)|mrm|/->`->`/[\Pi`\Pi S`R;s`r`]
\qtriangle(1500,-500)|mlm|/`<-`->/[P`\Pi S`\Pi;`p`sS]
\ptriangle(2000,-500)|mmr|/`->`<-/[\Pi S`R`\Pi;`sS`r]
\morphism(1700,100)|a|<75,-200>[`;\phi]
\morphism(2300,100)|a|<-75,-200>[`;\psi]
\efig
\end{equation}
Note that precomposing the top two triangles on the left above 
with $S$ produces the bottom two triangles on
the left above. These triangles give rise to pasting composites
that we have displayed, and named, on the right above. By
normality of $\bA$, and hence of $[\bA^2,\bA]$, we have 
$\phi\f p1_{\Pi}\ra p(sS)s$ and
$\psi\f r1_{\Pi}\ra r(sS)s$. It 
follows, by the fully faithful aspect of the universal property of 
the product cone of \ref{prt} that there is a unique 2-cell
\begin{equation}\label{thirteen}
\sigma\f1_{\Pi}\ra(sS)s
\end{equation}
with $p\sigma=\phi$ and 
$r\sigma=\psi$ in $[\bA^2,\bA]$. So $\sigma$ is by construction a
modification. Moreover, $\sigma$ is
invertible since $\phi$ and $\psi$ are so. If we precompose the
entire left side of (\ref{twelve}) with $S$ then we get a unique 2-cell
$\sigma S\f1_{\Pi S}\ra s(sS)$ with $(pS)(\sigma S)=\phi S$ and 
$(rS)(\sigma S)=\psi S$. It follows that $sS$ is an inverse 
equivalence of $s$.  We will take $s$ and the $\sigma$ 
as constructed here as the basic datum for a braiding (\ref{nine}) 
and the datum for a syllepsis (\ref{ten}).

\subsection{}\label{3.9} 
Before dealing with the further data for a braiding and the
braiding and syllepsis equations we will
establish the symmetry equation (\ref{eleven}) directly from the
descriptions of $\sigma$ and $\sigma S$. To show that 
$s\sigma=(\sigma S)s$ it suffices to show that 
$(rS)s\sigma=(rS)(\sigma S)s$ and $(pS)s\sigma=(pS)(\sigma S)s$.
We show the first equality by showing equality of the pasting 
composites of each side with the invertible $\mu\inv\f(rS)s\ra p$.
We have:
$$\bfig
\Vtriangle(0,0)/->`->`<-/<500,350>[\Pi`\Pi`\Pi S;1_{\Pi}`s`sS]
\morphism(1000,350)|a|<500,0>[\Pi`\Pi S;s]
\morphism(1500,350)|a|<500,0>[\Pi S`P;rS]
\morphism(1000,350)|b|/{@{->}@/_3em/}/<1000,0>[\Pi`P;p]
\morphism(500,250)|a|<0,-150>[`;\sigma]
\morphism(1500,250)|l|<0,-150>[`;\mu\inv]
\place(-150,-350)[=]
\morphism(0,-350)|a|<500,0>[\Pi`\Pi;1_{\Pi}]
\morphism(500,-350)|a|<500,0>[\Pi`\Pi S;s]
\morphism(1000,-350)|a|<500,0>[\Pi S`P;rS]
\morphism(0,-350)|l|<0,-500>[\Pi`\Pi S;s]
\morphism(500,-350)|b|/{@{->}@/_3em/}/<1000,0>[\Pi`P;p]
\morphism(1000,-450)|l|<0,-150>[`;\mu\inv]
\morphism(0,-850)|m|/{@{->}@/_3em/}/<1500,500>[\Pi S`P;rS]
\morphism(500,-650)|l|<0,-150>[`;\mu]
\morphism(0,-850)|l|<0,-500>[\Pi S`\Pi;sS]
\morphism(0,-1350)|b|/{@{->}@/_3em/}/<1500,1000>[\Pi`P;p]
\morphism(500,-1050)|l|<0,-150>[`;\nu S]
\place(-150,-1700)[=]
\morphism(0,-1700)|a|<500,0>[\Pi`\Pi S;s]
\Vtriangle(500,-2050)/->`->`<-/<500,350>[\Pi S`P`\Pi;rS`sS`p]
\morphism(1000,-1750)|a|<0,-150>[`;\nu S]
\place(-150,-2400)[=]
\morphism(0,-2400)|a|<500,0>[\Pi`\Pi S;s]
\morphism(500,-2400)|a|<500,0>[\Pi S`\Pi S;1_{\Pi S}]
\morphism(1000,-2400)|a|<1000,0>[\Pi S`P;rS]
\morphism(500,-2400)|l|<500,-350>[\Pi S`\Pi;sS]
\morphism(1000,-2750)|a|<1000,350>[\Pi`P;p]
\morphism(1000,-2750)|m|<500,0>[\Pi`\Pi S;s]
\morphism(1500,-2750)|m|<500,350>[\Pi S`P;rS]
\morphism(1000,-2500)|a|<0,-150>[`;\nu S]
\morphism(1500,-2600)|a|<0,-100>[`;\mu]
\morphism(1000,-2750)|b|/{@{->}@/_3em/}/<1000,350>[\Pi`P;p]
\morphism(1675,-2650)|r|<0,-100>[`;\mu\inv]
\place(-150,-3100)[=]
\morphism(0,-3100)|a|<500,0>[\Pi`\Pi S;s]
\Vtriangle(500,-3450)/->`->`<-/<500,350>%
[\Pi S`\Pi S`\Pi;1_{\Pi S}`sS`s]
\morphism(1500,-3100)|a|<500,0>[\Pi S`P;rS]
\morphism(1000,-3450)|b|/{@{->}@/_2em/}/<1000,350>[\Pi`P;p]
\morphism(1000,-3150)|l|<0,-150>[`;\sigma S]
\morphism(1500,-3250)|r|<0,-150>[`;\mu\inv]
\efig$$
The second equality, $(pS)s\sigma=(pS)(\sigma S)s$, can be shown, 
similarly, by showing equality of the pasting composites of each 
side with the invertible $\nu\inv\f(pS)s\ra r$.

\subsection{}\label{3.10} 
It is now convenient to again treat $X\x Y$, $Y\x X$, $X$, and $Y$ 
as objects in $[\bA^2,\bA]$, regarding them as alternative names for 
$\Pi=\x$, for $\Pi S$, and for the first and second projections of 
$\bA^2$ onto $\bA$. 

In $[\bA^3,\bA]$ we have (again writing $XY$ for $X\x Y$ and so on)
the composite
\begin{equation}
(XY)Z\to^a X(YZ)\to^s(YZ)X\to^a Y(ZX)\label{fourteen} 
\end{equation}
and the composite
\begin{equation}
(XY)Z\to^{sZ}(YX)Z\to^a Y(XZ)\to^{Ys}Y(ZX)\label{fifteen}
\end{equation}
We now write the composites in (\ref{fourteen}) and (\ref{fifteen}) 
as $m,n\f U\ra V$, where $U=\Pi(\Pi\x1)$ and 
$V=\Pi(1\x\Pi)(1\x S)(S\x 1)$, 
and we write $u\f U\ra\Pi_3$ and $v\f V\ra\Pi_3$ 
for the equivalences constructed using projections. 
Each of the arrows in (\ref{fourteen}) and (\ref{fifteen}) comes
with modifications constructed from either the $\mu$ in
(\ref{six}) or the $\mu$ and $\nu$ in (\ref{twelve}) and these fit together
and compose to give modifications $\alpha\f vm\ra u$ and
$\beta\f vn\ra u$. Now $(m,\alpha)$ and $(n,\beta)$ are arrows from
$U$ to $V$ in $[\bA^3,\bA](X,Y,Z)$, so by Proposition \ref{3.2} there
is a unique invertible 2-cell $(m,\alpha)\ra(n,\beta)$. This is the
modification `$R$' of \cite{McC}. The modification `$S$' of \cite{McC}
(which relates $s$ and $a^*$, the adjoint inverse equivalence of $a$)
is constructed in a similar way.

In \cite{McC} there are four diagrams, (BA1)--(BA4), to be satisfied 
in a braided monoidal category. Three of these demand
equality of two invertible 2-cells in $[\bA^4,\bA]$ (modifications)
and the other is an equality of two invertible 2-cells in 
$[\bA^3,\bA]$. In each case we use Proposition \ref{3.2} to
demonstrate equality. For the first three, the modifications in
question are shown to be 2-cells in $[\bA^4,\bA](X,Y,Z,W)$.
For the fourth, the modifications are shown to be 2-cells in 
$[\bA^3,\bA](X,Y,Z)$. 

In \cite{McC} there are two diagrams, (SA1) and (SA2), to be 
satisfied in a sylleptic (braided) monoidal category. Both
demand equality of two invertible 2-cells in $[\bA^3,\bA]$.
In each case, Proposition \ref{3.2} can be used after showing
that the modifications in question provide 2-cells in
$[\bA^3,\bA](X,Y,Z)$. As with the braiding equations, we forego
an explicit demonstration. To summarize:

\thm\label{finprodimpmonl}
A bicategory $\bA$ with binary product $\x$ and terminal object
$1$ underlies a symmetric monoidal bicategory with $\x$ as its
tensor product and $1$ as its unit object.
\frp
\eth

\section{Precartesian Bicategories}\label{precart}
\subsection{}\label{defprecart}
A bicategory $\bB$ is said to be {\em precartesian} if
\begin{enumerate}
\item[i)] the bicategory $\bM=\map\bB$ has finite products;
\item[ii)] each category $\bB(X,Y)$ has finite products.
\end{enumerate}
We henceforth assume $\bB$ to be a {\em precartesian} bicategory.

We have discussed the notations $\x\f\bM\x\bM\ra\bM$ and 
$p\f\x\ra P$ and $r\f\x\ra R$ for the binary product in a bicategory.
Here we find it convenient to write $I$ (rather than $1$) for the 
terminal object of $\bM$. For the binary product
in each $\bB(X,Y)$ we write $R\wedge S$, with $\pi\f R\wedge S\ra R$ 
and $\rho\f R\wedge S\ra S$ for the projections, and we use $\top$, or
$\top_{X,Y}$, for the terminal object of $\bB(X,Y)$. For the adjunction
units we write $\delta\f R\ra R\wedge R$ and $\tau\f R\ra\top$.

\exl\label{precarts}
\begin{enumerate}
\item[i)] All the examples of cartesian bicategories provided 
in \cite{caw};
\item[ii)] $\bB=\spn\E$, the bicategory of spans in $\E$, for $\E$  
a category with finite limits; 
\item[iii)] $\bB=\prof\E$, the bicategory of categories and profunctors
in $\E$, for $\E$ an elementary topos;
\item[iv)] $\bB=\cart\E$, the 2-category of categories with finite
products, finite product preserving functors, and natural 
transformations in $\E$, for $\E$ a category with finite limits.
\end{enumerate}
\eth

\subsection{}\label{precar} 
We now describe a certain bicategory $\bG$ associated to $\bB$. It is
in fact the bicategory formed by applying Street's \cite{str} two-sided 
Grothendieck construction to the pseudofunctor
\begin{equation}\label{Mhom}
\bM\op\x\bM\to^{i\op\x i}\bB\op\x\bB\to^{\bB(-,-)}\CAT
\end{equation}
wherein $i$ is the inclusion of $\bM$ in $\bB$ and $\bB(-,-)$ is
the hom pseudofunctor of (\ref{hom}). An object of $\bG$ is a triple
$(X,R,A)$ where $R$ is an object of the category $\bB(X,A)$ --- that 
is, a general arrow $R\f X\ra A$ in $\bB$. We shall often write
$R\f X\ra A$, or just $R$, for the object $(X,R,A)$. An arrow in $\bG$ 
from $(X,R,A)$ to $(Y,S,B)$ consists of a triple $(f,\alpha,u)$ where
$f\f X\ra Y$ and $u\f A\ra B$ are {\em maps} in $\bB$ and 
$\alpha\f uR\ra Sf$ is a 2-cell in $\bB$ as in
\begin{equation}\label{asquare}
\bfig
\square(0,0)[X`Y`A`B;f`R`S`u]
\morphism(125,250)|m|<250,0>[`;\alpha]
\efig
\end{equation}
and such arrows are composed by bicategorical pasting
as discussed in Verity's thesis \cite{ver}.

There is another way of describing an arrow of $\bG$.
Since the map $u$ has a right adjoint $u^*$, to give
the 2-cell $\alpha\f uR\ra Sf$ is equally to give a 2-cell
$\beta\f R\ra u^*Sf$ as in
\begin{equation}\label{bsquare}
\bfig
\square(0,0)/->`->`->`<-/[X`Y`A`B;f`R`S`u^*]
\morphism(125,250)|m|<250,0>[`;\beta]
\efig
\end{equation}
namely the {\em mate} of $\alpha$ in the sense of Kelly and 
Street \cite{kas}. We call the description of an
arrow of $\bG$ by a triple $(f,\alpha,u)$ as in (\ref{asquare})
its {\em primary form} and that by the triple $(f,\beta,u)$
as in (\ref{bsquare}) its {\em secondary form}. It is easy
to see that composition of arrows expressed in their secondary 
forms is again given by bicategorical pasting.

Given arrows $(f,\alpha,u)$ and $(f',\alpha',u')$ of $\bG$, expressed 
in primary form, a 2-cell in $\bG$ from $(f,\alpha,u)$ to 
$(f',\alpha',u')$ consists of 2-cells $\phi\f f\ra f'$ and 
$\psi\f u\ra u'$ in $\bM$ for which the diagram
\begin{equation}\label{a2-cell}
\bfig
\square(0,0)/->`<-`<-`->/[Sf`Sf'`uR`u'R;S\phi`\alpha`\alpha'`\psi R]
\efig
\end{equation}
commutes.
If the secondary forms of the same arrows are $(f,\beta, u)$ and
$(f',\beta',u')$ then the condition (\ref{a2-cell}) for 2-cells 
$\phi\f f\ra f'$ and $\psi\f u\ra u'$ 
to constitute a 2-cell in $\bG$ becomes the commutativity of
\begin{equation}\label{b2-cell}
\bfig
\Ctriangle(0,0)|amb|/<-``->/<350,350>[u^*Sf`R`u'^*Sf';\beta``\beta']
\Dtriangle(350,0)|mab|/`->`<-/<350,350>%
[u^*Sf`u^*Sf'`u'^*Sf';`u^*S\phi`\psi^*Sf']
\efig
\end{equation}
where $\psi^*\f u'^*\ra u^*$ is the mate with respect to the
adjunctions $u\laj u^*$ and $u'\laj u'^*$ of $\psi\f u\ra u'$.

\prp\label{Gequiv}
The typical arrow $(f,\alpha,u)$ of $G$, expressed in primary form,
is an equivalence if and only if $f$ and $u$ are equivalences
in $\bM$ and $\alpha$ is invertible in $\bB(X,B)$.
\frp
\eth

There are pseudofunctors
$\partial_0$ and $\partial_1$ as shown in (\ref{Gspan}) below
\begin{equation}\label{Gspan}
\bfig
\Atriangle/->`->`/[\bG`\bM`\bM;\partial_0`\partial_1`]
\efig
\end{equation}
given (using either description for arrows) by
$$\partial_0[(\phi,\psi)\f(f,\alpha,u)\ra(f',\alpha',u')%
\f(X,R,A)\ra(Y,S,B)]=\phi\f f\ra f'\f X\ra Y$$
and 
$$\partial_1[(\phi,\psi)\f(f,\alpha,u)\ra(f',\alpha',u')%
\f(X,R,A)\ra(Y,S,B)]=\psi\f u\ra u'\f A\ra B$$
We note that $\partial_0$ and $\partial_1$ have
the common section $\iota\f\bM\ra\bG$ which sends a map
$f\f X\ra Y$ to the square
$$\bfig
\square(0,0)[X`Y`X`Y;f`1_X`1_Y`f]
\efig$$
with the identity 2-cell understood (and meaningful by the normality 
of $\bB$).

\subsection{}\label{xinG} 
We are going to show that the bicategory $\bG$ has finite products. 
We shall use $R\ox S$ for the binary product in $\bG$, to distinguish
$f\ox g$ from $f\x g$ when $f$ and $g$ are maps.

Given objects $R\f X\ra A$ and $S\f Y\ra B$ of $\bG$, we define
$R\ox S\f X\x Y\ra A\x B$ by
$$R\ox S = p^*_{A,B}Rp_{X,Y}\wedge r^*_{A,B}Sr_{X,Y}$$
which can be abbreviated to $p^*Rp\wedge r^*Sr$ when $X$, $Y$, $A$,
and $B$ are clear; and we define arrows $p_{R,S}\f R\ox S\ra R$
and $r_{R,S}\f R\ox S\ra S$ of $\bG$, in their secondary forms, by
$$p_{R,S}=(p_{X,Y},\pi,p_{A,B})\quad\mbox{and}\quad%
r_{R,S}=(r_{X,Y},\rho,r_{A,B})$$
as in the diagram
\begin{equation}\label{projexionsb}
\bfig
\square(500,0)|almb|/<-`->`->`->/%
[X`X\x Y`A`A\x B;p_{X,Y}`R`R\ox S`p^*_{A,B}]
\morphism(850,250)|m|<-200,0>[`;\pi]
\square(1000,0)|amrb|/->`->`->`<-/%
[X\x Y`Y`A\x B`B;r_{X,Y}`R\ox S`S`r^*_{A,B}]
\morphism(1150,255)|m|<200,0>[`;\rho]
\efig
\end{equation}
wherein $\pi$ and $\rho$ are the projections of the product
$\wedge$ in $\bB(X\x Y,A\x B)$. Accordingly, the primary forms of 
$p_{R,S}$ and of $r_{R,S}$ are $(p_{X,Y},\tilde{p}_{R,S},p_{A,B})$
and $(r_{X,Y},\tilde{r}_{R,S},r_{A,B})$ as in
\begin{equation}\label{projexionsa}
\bfig
\square(500,0)|almb|/<-`->`->`<-/%
[X`X\x Y`A`A\x B;p_{X,Y}`R`R\ox S`p_{A,B}]
\morphism(850,250)|m|<-250,0>[`;\tilde{p}_{R,S}]
\square(1000,0)|amrb|[X\x Y`Y`A\x B`B;r_{X,Y}`R\ox S`S`r_{A,B}]
\morphism(1150,255)|m|<250,0>[`;\tilde{r}_{R,S}]
\efig
\end{equation}
where the 2-cell $\tilde{p}_{R,S}$ is the mate 
$p(p^*Rp\wedge r^*Sr){\s\to^{p\pi}}pp^*Rp{\s\to^{\epsilon_pRp}}Rp$ 
of $\pi$ and similarly for $\tilde{r}_{R,S}$. 

To prove that $R\ox S$, with the projections $p_{R,S}$ and
$r_{R,S}$, is indeed the binary product in $\bG$ is to show
that, for any object $T\f Z\ra C$ of $\bG$, the functor
$$\bG(T,R\ox S)\to^{(\bG(T,p_{R,S}),\bG(T,r_{R,S}))}%
\bG(T,R)\x\bG(T,S)$$
is essentially surjective on objects and is fully faithful. These
properties are established in the following lemmas, in the first
of which we use the fact that a functor of the form $\bB(f,u^*)$,
where $f$ and $u$ are maps, preserves products, being the right adjoint
of $\bB(f^*,u)$.

\lem\label{eso}
For each object $T\f Z\ra C$ in $\bG$, the functor
$$\bG(T,R\ox S)\to^{(\bG(T,p_{R,S}),\bG(T,r_{R,S}))}%
\bG(T,R)\x\bG(T,S)$$
is essentially surjective on objects. In fact, more is true:
for any 
$$(f,\alpha,u)\f(Z,T,C)\ra(X,R,A)\quad\mbox{and}\quad 
(g,\beta,v)\f(Z,T,C)\ra(Y,S,B)$$ in $\bG$ and any
$$(h,w)\f(Z,C)\ra(X\x Y,A\x B)$$ with invertible 2-cells
$$(\mu_0,\mu_1)\f(p_{X,Y},p_{A,B})(h,w)\ra(f,u)\quad\mbox{and}\quad
(\nu_0,\nu_1)\f(r_{X,Y},r_{A,B})(h,w)\ra(g,v)$$ provided by the products
in $\bM\x\bM$, there
exists a unique $\gamma$ making $(h,\gamma,w)$ an arrow
$T\ra R\ox S$ in $\bG$ with invertible 2-cells
$$(\mu_0,\mu_1)\f p_{R,S}(h,\gamma,w)\ra(f,\alpha,u)%
\quad\mbox{and}\quad%
(\nu_0,\nu_1)\f r_{R,S}(h,\gamma,w)\ra(g,\beta,v).$$
\eth
\prf
Here we use the secondary forms for all arrows in $\bG$. 
The product diagram
$$p^*Rp\toleft^{\pi}p^*Rp\wedge r^*Sr\to^{\rho}r^*Sr$$ 
in $\bB(X\x Y,A\x B)$ gives rise to the further product diagram
$${\scalefactor{0.75}w^*p^*Rf\toleft^{w^*p^*R\mu_0}w^*p^*Rph%
\toleft^{w^*\pi h}w^*(p^*Rp\wedge r^*Sr)h\to^{w^*\rho h}w^*r^*Srh%
\to^{w^*r^*S\nu_0}w^*r^*Sg}$$
in $\bB(Z,C)$
because $\bB(h,w^*)$ preserves products
and $\mu_0$ and $\nu_0$ are invertible. Now consider the following
diagram: 
$$\bfig
\square(0,0)/<-`<-``/<750,500>%
[w^*p^*Rf`w^*p^*Rph`u^*Rf`;w^*p^*R\mu_0`\mu_1^*Rf``]
\square(750,0)|almb|/<-``<--`/<750,500>%
[w^*p^*Rph`w^*(R\ox S)h``T;w^*\pi h``\gamma`]
\morphism(1500,0)|b|<-1500,0>[T`u^*Rf;\alpha]
\square(1500,0)|amrb|/->`<--``/<750,500>%
[w^*(R\ox S)h`w^*r^*Srh`T`;w^*\rho h`\gamma``]
\square(2250,0)/->``<-`/<750,500>%
[w^*r^*Srh`w^*r^*Sg``v^*Sg;w^*r^*S\nu_0``\nu^*_1Sg`]
\morphism(1500,0)|b|<1500,0>[T`v^*Sg;\beta]
\efig$$
To give such a 2-cell $\gamma$ is to give an arrow $(h,\gamma,w)\f T\ra R\ox S$
in $\bG$; and by (\ref{b2-cell}), the condition for
$(\mu_0,\mu_1)\f p_{R,S}(h,\gamma,w)\ra(f,\alpha,u)$ to
provide a~2-cell in $\bG$ is the commutativity of the
left square, while the condition for
$(\nu_0,\nu_1)\f r_{R,S}(h,\gamma,w)\ra(g,\beta,v)$ to
provide a \linebreak 2-cell in $\bG$ is the commutativity of the
right square. Since the top row of the diagram is a product,
there is a unique such $\gamma$.
\frp

\cor\label{!ness}
If two arrows 
$(h,\gamma,w),(h,\delta,w)\f(Z,T,C)\ra(X\x Y,R\ox S, A\x B)$
in $\bG$ satisfy $p_{R,S}(h,\gamma,w)=p_{R,S}(h,\delta,w)$
and $r_{R,S}(h,\gamma,w)=r_{R,S}(h,\delta,w)$ then 
$\gamma=\delta$ so that $(h,\gamma,w)=(h,\delta,w)$.
\eth
\prf
We can apply Lemma \ref{eso} with $(f,\alpha,u)=p_{R,S}(h,\gamma,w)$ 
and with $(g,\beta,v)=r_{R,S}(h,\gamma,w)$, taking $\mu_0$, $\mu_1$, 
$\nu_0$, and $\nu_1$ to be identities. 
\frp

\lem\label{f&f}
For each object $T\f Z\ra C$ in $\bG$, the functor
$$\bG(T,R\ox S)\to^{(\bG(T,p_{R,S}),\bG(T,r_{R,S}))}%
\bG(T,R)\x\bG(T,S)$$
is fully faithful.
\eth
\prf
Let $(h,\gamma,w),(k,\delta,x)\f(Z,T,C)\ra(X\x Y,R\ox S, A\x B)$ 
be arrows of $\bG$ in primary form, and consider 2-cells 
$$(\phi,\psi)\f p_{R,S}(h,\gamma,w)\ra p_{R,S}(k,\delta,x)%
\quad\mbox{and}\quad%
(\chi,\omega)\f r_{R,S}(h,\gamma,w)\ra r_{R,S}(k,\delta,x)$$
Since these data further provide 2-cells
$(\phi,\psi)\f (p_{X,Y},p_{A,B})(h,w)\ra (p_{X,Y},p_{A,B})(k,x)$ and
$(\chi,\omega)\f (r_{X,Y},r_{A,B})(h,w)\ra (r_{X,Y},r_{A,B})(k,x)$
in $\bM\x\bM$ and since the bicategory $\bM\x\bM$ has finite products, 
there are unique 2-cells $\langle\phi,\chi\rangle\f h\ra k$
and $\langle\psi,\omega\rangle\f w\ra x$ in $\bM$ satisfying
$$(p_{X,Y},p_{A,B})(\langle\phi,\chi\rangle,\langle\psi,\omega\rangle)%
=(\phi,\psi)\quad\mbox{and}\quad%
(r_{X,Y},r_{A,B})(\langle\phi,\chi\rangle,\langle\psi,\omega\rangle)%
=(\chi,\omega).$$
It only remains to show that 
$(\langle\phi,\chi\rangle,\langle\psi,\omega\rangle)$ constitutes a 
2-cell $(h,\gamma,w)\ra(k,\delta,x)$ in~$\bG$, for which the requisite 
condition is
\begin{equation}\label{=pr}
\bfig
\square(0,0)/{@{->}@/^1em/}`->`->`{@{->}@/_1em/}/%
[Z`X\x Y`C`A\x B;k`T`R\ox S`w]
\morphism(0,500)|b|/{@{->}@/_1em/}/<500,0>[Z`X\x Y;h]
\morphism(300,450)|l|/->/<0,125>[`;\langle\phi,\chi\rangle]
\morphism(175,150)|a|/->/<150,0>[`;\gamma]
\place(1000,250)[=]
\square(1500,0)/{@{->}@/^1em/}`->`->`{@{->}@/_1em/}/%
[Z`X\x Y`C`A\x B;k`T`R\ox S`w]
\morphism(1500,0)|a|/{@{->}@/^1em/}/<500,0>[C`A\x B;x]
\morphism(1675,325)|a|/->/<150,0>[`;\delta]
\morphism(1800,-75)|l|/->/<0,150>[`;\langle\psi,\omega\rangle]
\efig
\end{equation}
Because both pasting composites displayed are of the form
$(k,?,w)\f T\ra R\ox S$ in $\bG$, Corollary~\ref{!ness} ensures
their equality if they are coequalized by both $p_{R,S}$ 
and $r_{R,S}$.
Composing the left side of (\ref{=pr}) with $p_{R,S}$ gives the
left pasting composite below, while composing the right side of
(\ref{=pr}) with $p_{R,S}$ gives the right pasting composite below.
$$\bfig
\square(0,0)|almb|[Z`X\x Y`C`A\x B;h`T`R\ox S`w]
\morphism(175,250)|a|/->/<150,0>[`;\gamma]
\square(500,0)|amrb|[X\x Y`X`A\x B`A;p_{X,Y}`R\ox S`R`p_{A,B}]
\morphism(675,250)|a|/->/<150,0>[`;\tilde{p}_{R,S}]
\square(2000,0)|almb|[Z`X\x Y`C`A\x B;k`T`R\ox S`x]
\morphism(2175,250)|a|/->/<150,0>[`;\delta]
\square(2500,0)|amrb|[X\x Y`X`A\x B`A;p_{X,Y}`R\ox S`R`p_{A,B}]
\morphism(2675,250)|a|/->/<150,0>[`;\tilde{p}_{R,S}]
\morphism(0,500)|a|/{@{->}@/^1em/}/<500,250>[Z`X\x Y;k]
\morphism(500,750)|a|/{@{->}@/^1em/}/<500,-250>[X\x Y`X;p_{X,Y}]
\morphism(500,575)|l|/->/<0,100>[`;\phi]
\morphism(2000,0)|b|/{@{->}@/_1em/}/<500,-250>[C`A\x B;w]
\morphism(2500,-250)|b|/{@{->}@/_1em/}/<500,250>[A\x B`A;p_{A,B}]
\morphism(2500,-175)|l|/->/<0,100>[`;\psi]
\efig$$
However, the composites above are equal since $(\phi,\psi)$ was
assumed to provide a 2-cell in~$\bG$. A similar
argument applies to composition with $r_{R,S}$ and this completes
the proof.
\frp
\thm\label{Ghasx}
$\bG$ has finite products preserved by $\partial_0$ and $\partial_1$. 
\eth
\prf
The binary products are provided by Lemmas \ref{eso} and \ref{f&f};
we now show that $\top=\top_{I,I}\f I\ra I$ provides a terminal object 
for $\bG$.  In fact, for each object $(X,R,A)$ in 
$\bG$ we have the maps $t_X\f X\ra I$ and $t_A\f A\ra I$ of \ref{prt};
and the right adjoint functor $\bB(t_X,t_A^*)$ sends $\top_{I,I}$ to 
$\top_{X,A}$. Accordingly, there is a
unique 2-cell $\tau\f R\ra t_A^*\top t_X$, which is  the
description in secondary form of an arrow $t_R\f R\ra \top$ in $\bG$. 
For any arrow
$(f,\alpha,u)\f (X,R,A)\ra (Y,S,B)$ we have a unique 2-cell
$t_S(f,\alpha,u)\ra t_R$ given by $(t'_f,t'_u)$ and it is
invertible.
\frp

\rmk\label{tunit}
For the primary form of the arrow $t_R$ we use $(t_X,\tilde{t}_R,t_A)$.
It is the component of a pseudonatural transformation 
$t\f1_{\bG}\ra\top!$
which is the unit for a pseudoadjunction $!\laj\top\f\one\ra\bG$.
\eth

\subsection{}\label{dunit}
When $Y=X$ and $B=A$ in (\ref{projexionsa}), we have 
$d^*_A(R\ox S)d_X=d^*_A(p^*Rp\wedge r^*Sr)d_X%
\cong d^*_Ap^*Rpd_X\wedge d^*_Ar^*Srd_X\cong R\wedge S$; thus 
$$R\wedge S\cong d^*_A(R\ox S)d_X$$
We have in particular $R\wedge R \iso d^*_A(R\ox R)d_X$; and 
composing this isomorphism with $\delta\f R\ra R\wedge R$ gives 
a 2-cell we can still call $\delta$ as in
$$\bfig
\square(0,0)/->`->`->`<-/[X`X\x X`A`A\x A;d_X`R`R\ox R`d^*_A]
\morphism(125,250)|m|<250,0>[`;\delta]
\efig$$
giving, in secondary form, an arrow $d_R\f R\ra R\ox R$, for
whose primary form we write $(d_X,\tilde{d}_R,d_A)$. The $d_R$ are
the components of a pseudonatural
transformation $d\f1_{\bG}\ra\ox\Delta$, 
which is the unit for a pseudoadjunction 
$\Delta\laj\ox:\bG\x\bG\ra\bG$ (for which the counit has
components given by the pairs $(p_{R,S},r_{R,S})=%
((p_{X,Y},\tilde p_{R,S},p_{A,B}),(r_{X,Y},\tilde r_{R,S},r_{A,B}))$ 
of (\ref{projexionsa})). We stress that the $\tilde p_{R,S}$,
$\tilde r_{R,S}$, $\tilde d_R$, and $\tilde t_R$ are not,
in general, invertible.

\subsection{}\label{aoxb}
Each 2-cell $\alpha\f R\ra R'\f X\ra A$ in $\bB$ gives an arrow
$(1_X,\alpha,1_A)$ of $\bG$, as in
$$\bfig
\square(0,0)[X`X`A`A;1`R`R'`1]
\morphism(125,250)|m|<250,0>[`;\alpha]
\efig$$
If $\gamma\f R'\ra R''\f X\ra A$ is another 2-cell, the composite
arrow $(1,\gamma,1)(1,\alpha,1)$ is $(1,\gamma\alpha,1)$, since
the normality of $\bB$ forces the constraint $(1R')1\iso 1(R'1)$
to be an identity.

Like $\x$ in \ref{prt} the product $\ox$ of $\bG$ provides a 
pseudofunctor $\ox\f\bG\x\bG\ra\bG$. Consider the product
$(1_X,\alpha,1_A)\ox(1_Y,\beta,1_B)$, where $\alpha\f R\ra R'\f X\ra A$
and $\beta\f S\ra S'\f Y\ra B$. Like the $f\x g$ of \ref{prt},
it is determined to within isomorphism by the existence of
invertible 2-cells
$$\mu: p_{R',S'}(1_X,\alpha,1_A)\ox(1_Y,\beta,1_B)\iso%
(1_X,\alpha,1_A)p_{R,S},$$
$$\nu: r_{R',S'}(1_X,\alpha,1_A)\ox(1_Y,\beta,1_B)\iso%
(1_Y,\beta,1_B)r_{R,S}.$$
Recall from \ref{xinG} that, in secondary form, 
$p_{R,S}=(p_{X,Y},\pi,p_{A,B})$ and 
$r_{R,S}=(r_{X,Y},\rho,r_{A,B})$. Since we have equalities 
$$(p_{X,Y},p_{A,B})(1_{X\x Y},1_{A\x B})=(1_X,1_A)(p_{X,Y},p_{A,B}),$$
$$(r_{X,Y},r_{A,B})(1_{X\x Y},1_{A\x B})=(1_Y,1_B)(r_{X,Y},r_{A,B}),$$
it follows from Lemma \ref{eso} that there is a unique
$\phi\f R\ox S\ra R'\ox S'$ for which we have equalities
\begin{equation}\label{lnat1}
p_{R',S'}(1_{X\x Y},\phi,1_{A\x B})=%
(1_X,\alpha,1_A)p_{R,S}
\end{equation}
\begin{equation}\label{lnat2}
r_{R',S'}(1_{X\x Y},\phi,1_{A\x B})=%
(1_Y,\beta,1_B)r_{R,S}
\end{equation}
We write $\alpha\ox\beta$ for this value of $\phi$.
Inserting the values above of $\tilde{p}_{R,S}$ and so on in
these last equalities gives 
$$\pi(\alpha\ox\beta)=(p^*\alpha p)\pi$$
$$\rho(\alpha\ox\beta)=(r^*\beta r)\rho$$
from which we deduce that
$$\alpha\ox\beta=(p^*\alpha p)\wedge(r^*\beta r)$$
Thus our formula $R\ox S=(p^*R p)\wedge(r^*S r)$ extends
to 2-cells to give a functor 
$$\ox_{(X,Y),(A,B)}\f \bB\x\bB((X,Y),(A,B))=%
\bB(X,A)\x\bB(Y,B)\ra\bB(X\x Y,A\x B)$$
namely the composite 
$$\bB(X,A)\x\bB(Y,B){\s\to^{\bB(p,p^*)\x\bB(r,r^*)}}%
\bB(X\x Y,A\x B)\x\bB(X\x Y,A\x B)%
{\s\to^{\wedge}}\bB(X\x Y,A\x B)$$

\rmk\label{lxnat}
We observe for later reference that (\ref{lnat1}) and (\ref{lnat2}) can be displayed
as
$$\bfig
\square(0,0)/->```->/[X\ox Y`X`A\ox B`A;p_{X,Y}```p_{A,B}]
\morphism(0,500)|l|/{@{->}@/_2em/}/<0,-500>[X\ox Y`A\ox B;R\ox S]
\morphism(0,500)|m|/{@{->}@/^2em/}/<0,-500>[X\ox Y`A\ox B;R'\ox S']
\morphism(500,500)|r|/{@{->}@/^2em/}/<0,-500>[X`A;R']
\morphism(-150,250)|a|<200,0>[`;\alpha\ox\beta]
\morphism(350,250)|a|<200,0>[`;\tilde p_{R',S'}]
\place(1000,500)[=]
\square(1500,0)/->```->/[X\ox Y`X`A\ox B`A;p_{X,Y}```p_{A,B}]
\morphism(1500,500)|l|/{@{->}@/_2em/}/<0,-500>[X\ox Y`A\ox B;R\ox S]
\morphism(2000,500)|l|/{@{->}@/_2em/}/<0,-500>[X`A;R]
\morphism(2000,500)|r|/{@{->}@/^2em/}/<0,-500>[X`A;R']
\morphism(1450,250)|a|<200,0>[`;\tilde p_{R,S}]
\morphism(1900,250)|a|<200,0>[`;\alpha]
\efig$$
and
$$\bfig
\square(0,0)/->```->/[X\ox Y`Y`A\ox B`B;r_{X,Y}```r_{A,B}]
\morphism(0,500)|l|/{@{->}@/_2em/}/<0,-500>[X\ox Y`A\ox B;R\ox S]
\morphism(0,500)|m|/{@{->}@/^2em/}/<0,-500>[X\ox Y`A\ox B;R'\ox S']
\morphism(500,500)|r|/{@{->}@/^2em/}/<0,-500>[Y`B;S']
\morphism(-150,250)|a|<200,0>[`;\alpha\ox\beta]
\morphism(350,250)|a|<200,0>[`;\tilde r_{R',S'}]
\place(1000,500)[=]
\square(1500,0)/->```->/[X\ox Y`Y`A\ox B`B;r_{X,Y}```r_{A,B}]
\morphism(1500,500)|l|/{@{->}@/_2em/}/<0,-500>[X\ox Y`A\ox B;R\ox S]
\morphism(2000,500)|l|/{@{->}@/_2em/}/<0,-500>[Y`B;S]
\morphism(2000,500)|r|/{@{->}@/^2em/}/<0,-500>[Y`B;S']
\morphism(1450,250)|a|<200,0>[`;\tilde r_{R,S}]
\morphism(1900,250)|a|<200,0>[`;\beta]
\efig$$
saying that the $\tilde p_{R,S}$ and $\tilde r_{R,S}$ are natural
in $R$ and $S$.
\eth

\subsection{}\label{lax}
We extend the definition of $\ox$ to objects by setting
$X\ox Y= X\x Y$.
We shall now show that the functors 
$$\ox_{(X,Y),(A,B)}\f\bB(X,A)\x\bB(Y,B)\ra\bB(X\ox Y,A\ox B)$$
provide the effect on homs for a {\em lax} functor
$$(\ox,\widetilde\ox,\ox^\circ)\f\bB\x\bB\ra\bB$$ 

First, for objects $X$ and $Y$ of $\bB$, there is by
Lemma \ref{eso} a unique $\ox^\circ\f 1_{X\ox Y}\ra 1_X\ox 1_Y$
satisfying
$$\bfig
\square(0,0)|almb|[X\ox Y`X\ox Y`X\ox Y`X\ox Y;%
1`1_{X\ox Y}`1_X\ox 1_Y`1]
\morphism(125,250)|m|<200,0>[`;\ox^\circ]
\square(500,0)|amrb|[X\ox Y`X`X\ox Y`X;p`1_X\ox1_Y`1_X`p]
\morphism(650,250)|m|<200,0>[`;\tilde{p}_{1,1}]
\place(1500,500)[=]
\square(2000,0)[X\ox Y`X`X\ox Y`X;p`1_{X\ox Y}`1_X`p]
\efig$$
and

$$\bfig
\square(0,0)|almb|[X\ox Y`X\ox Y`X\ox Y`X\ox Y;%
1`1_{X\ox Y}`1_X\ox 1_Y`1]
\morphism(125,250)|m|<200,0>[`;\ox^\circ]
\square(500,0)|amrb|[X\ox Y`Y`X\ox Y`Y;r`1_X\ox1_Y`1_Y`r]
\morphism(650,250)|m|<200,0>[`;\tilde{r}_{1,1}]
\place(1500,500)[=]
\square(2000,0)[X\ox Y`Y`X\ox Y`Y;r`1_{X\ox Y}`1_Y`r]
\efig$$

Given $R\f X\ra A$, $S\f Y\ra B$, 
$T\f A\ra L$ and $U\f B\ra M$, vertical pasting as in the diagram
\begin{equation}\label{tentilde}
\bfig
\square(500,500)|almm|/<-`->`->`<-/[X`X\ox Y`A`A\ox B;p`R`R\ox S`p]
\morphism(850,750)|a|<-200,0>[`;\tilde{p}_{R,S}]
\square(500,0)|mlmb|/<-`->`->`<-/[A`A\ox B`L`L\ox M;p`T`T\ox U`p]
\morphism(850,250)|a|<-200,0>[`;\tilde{p}_{T,U}]
\square(1000,500)|amrm|[X\ox Y`Y`A\ox B`B;r`R\ox S`S`r]
\morphism(1150,750)|a|<200,0>[`;\tilde{r}_{R,S}]
\square(1000,0)|mmrb|[A\ox B`B`L\ox M`M;r`T\ox U`U`r]
\morphism(1150,250)|a|<200,0>[`;\tilde{r}_{T,U}]
\efig
\end{equation}
gives arrows in $\bG$ from $(T\ox U)( R\ox S)$ to $TR$ and $US$. 
Accordingly, there is by Lemma~\ref{eso} a unique
$$\widetilde\ox : (T\ox U)(R\ox S)\ra (TR)\ox(US)$$
whose composites with $p_{TR,US}$ and $r_{TR,US}$ are these vertical
pastings.

The first requirement for $(\ox,\widetilde\ox,\ox^\circ)$ to be a 
lax functor $\bB\x\bB\ra\bB$ is that $\widetilde\ox$ be natural in 
$T$, $R$, $U$, and $S$; and it is so because the assignment 
$\alpha\rass(1,\alpha,1)$ respects {\em both} vertical {\em and}
horizontal composition of 2-cells. 

For the associativity coherence condition, consider the data in 
(\ref{tentilde}) along with further arrows $V\f L\ra C$ and 
$W\f M\ra D$ in $\bB$. We require the following equality of pasting
composites:

\begin{equation}\bfig\label{assoc}
\square(0,1000)|almm|/->`->``/[\bullet`\bullet`\bullet`;%
1`R\ox S``]
\square(0,500)|almm|/`->``/[\bullet``\bullet`\bullet;%
`T\ox U``]
\morphism(500,1500)|m|<0,-1000>[\bullet`\bullet;TR\ox US]
\square(0,0)|mlmb|/->`->`->`->/[\bullet`\bullet`\bullet`\bullet;%
1`V\ox W`V\ox W`1]
\morphism(125,1000)|a|<150,0>[`;\widetilde\ox]
\morphism(500,1500)|a|<500,0>[\bullet`\bullet;1]
\morphism(500,0)|b|<500,0>[\bullet`\bullet;1]
\morphism(1000,1500)|r|<0,-1500>[\bullet`\bullet;VTR\ox WUS]
\morphism(675,750)|a|<150,0>[`;\widetilde\ox]
\place(1500,1500)[=]
\square(2000,1000)|almm|/->`->`->`->/%
[\bullet`\bullet`\bullet`\bullet;1`R\ox S`R\ox S`1]
\square(2000,500)|almm|/`->``/[\bullet`\bullet`\bullet`;%
`T\ox U``]
\morphism(2500,1000)|m|<0,-1000>[\bullet`\bullet;VT\ox WU]
\square(2000,0)|mlmb|/`->``->/[\bullet``\bullet`\bullet;%
`V\ox W``1]
\morphism(2125,500)|a|<150,0>[`;\widetilde\ox]
\morphism(2500,1500)|a|<500,0>[\bullet`\bullet;1]
\morphism(2500,0)|b|<500,0>[\bullet`\bullet;1]
\morphism(3000,1500)|r|<0,-1500>[\bullet`\bullet;VTR\ox WUS]
\morphism(2675,750)|a|<150,0>[`;\widetilde\ox]
\efig\end{equation}

To establish this equality consider the effect of pasting
$p_{VTR,WUS}$ on the right to the left side of Equation~(\ref{assoc}):
$$\bfig
\scalefactor{0.93}
\square(0,1000)|ammm|/->`->``/[\bullet`\bullet`\bullet`;%
1`R\ox S``]
\square(0,500)|ammm|/`->``/[\bullet``\bullet`\bullet;%
`T\ox U``]
\morphism(500,1500)|m|<0,-1000>[\bullet`\bullet;TR\ox US]
\square(0,0)|mmmb|/->`->`->`->/[\bullet`\bullet`\bullet`\bullet;%
1`V\ox W`V\ox W`1]
\morphism(125,1000)|a|<150,0>[`;\widetilde\ox]
\morphism(500,1500)|a|<500,0>[\bullet`\bullet;1]
\morphism(500,0)|b|<500,0>[\bullet`\bullet;1]
\morphism(1000,1500)|m|<0,-1500>[\bullet`\bullet;VTR\ox WUS]
\morphism(625,750)|a|<150,0>[`;\widetilde\ox]
\square(1000,0)|amrb|/->`->`->`->/<500,1500>%
[\bullet`\bullet`\bullet`\bullet;
p`VTR\ox WUS`VTR`p]
\morphism(1150,500)|a|<200,0>[`;\tilde{p}_{VTR,WUS}]
\place(1750,1500)[=]
\square(2000,1000)|ammm|/->`->``/[\bullet`\bullet`\bullet`;%
1`R\ox S``]
\square(2000,500)|ammm|/`->``/[\bullet``\bullet`\bullet;%
`T\ox U``]
\morphism(2500,1500)|m|<0,-1000>[\bullet`\bullet;TR\ox US]
\square(2000,0)|mmmb|/->`->`->`->/[\bullet`\bullet`\bullet`\bullet;%
1`V\ox W`V\ox W`1]
\morphism(2100,1000)|a|<150,0>[`;\widetilde\ox]
\square(2500,500)|amrm|/->`->`->`->/<500,1000>[\bullet`\bullet`\bullet`\bullet;%
p`TR\ox US`TR`p]
\square(2500,0)|mmrb|/->`->`->`_>/[\bullet`\bullet`\bullet`\bullet;%
p`V\ox W`V`p]
\morphism(2625,750)|a|<250,0>[`;\tilde{p}_{TR,US}]
\morphism(2650,250)|a|<200,0>[`;\tilde{p}_{V,W}]
\place(3250,1500)[=]
\square(3500,1000)|almm|/->`->`->`->/%
[\bullet`\bullet`\bullet`\bullet;p`R\ox S`R`p]
\square(3500,500)|mlmm|/->`->`->`->/%
[\bullet`\bullet`\bullet`\bullet;p`T\ox U`T`p]
\square(3500,0)|mlmb|/->`->`->`->/%
[\bullet`\bullet`\bullet`\bullet;p`V\ox W`V`p]
\morphism(3650,1250)|a|<200,0>[`;\tilde{p}_{R,S}]
\morphism(3650,750)|a|<200,0>[`;\tilde{p}_{T,U}]
\morphism(3650,250)|a|<200,0>[`;\tilde{p}_{V,W}]
\efig$$
It is clear that when $p_{VTR,WUS}$ is pasted to the right side
of Equation~(\ref{assoc}) the result is the same. Similarly, the left
side of (\ref{assoc}) pasted to $r_{VTR,WUS}$ is equal to
the right side of (\ref{assoc}) pasted to $r_{VTR,WUS}$. 
This suffices by Corollary \ref{!ness} to prove the condition 
satisfied.  

For the unitary coherence conditions we require:

\begin{equation}\bfig\label{unitary}
\scalefactor{0.97}
\square(0,500)|almm|/->`->`->`->/[X\ox Y`X\ox Y`X\ox Y`X\ox Y;%
1`1_{X\ox Y}`1_X\ox1_Y`1]
\square(0,0)|mlmb|/->`->`->`->/[X\ox Y`X\ox Y`A\ox B`A\ox B;%
1`R\ox S`R\ox S`1]
\square(500,0)|alrb|/->``->`->/<500,1000>%
[X\ox Y`X\ox Y`A\ox B`A\ox B;1``R\ox S`1]
\morphism(125,750)|a|<150,0>[`;\ox^\circ]
\morphism(700,500)|a|<150,0>[`;\widetilde\ox]
\place(1250,1000)[=]
\square(1500,0)|ammb|/->`->`->`->/<500,1000>%
[X\ox Y`X\ox Y`A\ox B`A\ox B;1`R\ox S`R\ox S`1]
\place(2250,1000)[=]
\square(2500,500)|almm|/->`->`->`->/[X\ox Y`X\ox Y`A\ox B`A\ox B;%
1`R\ox S`R\ox S`1]
\square(2500,0)|mlmb|/->`->`->`->/[A\ox B`A\ox B`A\ox B`A\ox B;%
1`1_{A\ox B}`1_A\ox1_B`1]
\square(3000,0)|alrb|/->``->`->/<500,1000>%
[X\ox Y`X\ox Y`A\ox B`A\ox B;1``R\ox S`1]
\morphism(2625,250)|a|<150,0>[`;\ox^\circ]
\morphism(3200,500)|a|<150,0>[`;\widetilde\ox]
\efig\end{equation}

To prove the first of these consider
$$\bfig
\scalefactor{0.70}
\square(0,500)|almm|/->`->`->`->/[\bullet`\bullet`\bullet`\bullet;%
1`1_{X\ox Y}`1_X\ox1_Y`1]
\square(0,0)|mlmb|/->`->`->`->/[\bullet`\bullet`\bullet`\bullet;%
1`R\ox S`R\ox S`1]
\square(500,0)|almb|/->``->`->/<500,1000>%
[\bullet`\bullet`\bullet`\bullet;1``R\ox S`1]
\morphism(125,750)|a|<150,0>[`;\ox^\circ]
\morphism(700,500)|a|<150,0>[`;\widetilde\ox]
\square(1000,0)|amrb|/->``->`->/<500,1000>%
[\bullet`\bullet`\bullet`\bullet;p``R`p]
\morphism(1200,500)|a|<150,0>[`;\tilde{p}_{R,S}]
\place(1750,1000)[=]
\square(2000,500)|almm|/->`->`->`->/%
[\bullet`\bullet`\bullet`\bullet;1`1_{X\ox Y}`1_X\ox1_Y`1]
\square(2000,0)|mlmb|/->`->`->`->/%
[\bullet`\bullet`\bullet`\bullet;1`R\ox S`R\ox S`1]
\morphism(2125,750)|a|<150,0>[`;\ox^\circ]
\square(2500,0)|mmrb|/->``->`->/%
[\bullet`\bullet`\bullet`\bullet;p``R`p]
\morphism(2700,250)|a|<150,0>[`;\tilde{p}_{R,S}]
\square(2500,500)|amrm|/->``->`->/%
[\bullet`\bullet`\bullet`\bullet;p``1_X`p]
\morphism(2700,750)|a|<150,0>[`;\tilde{p}_{1_X,1_Y}]
\place(3250,1000)[=]
\square(3500,0)|mlrb|/->`->`->`->/%
[\bullet`\bullet`\bullet`\bullet;p`R\ox S`R`p]
\morphism(3700,250)|a|<150,0>[`;\tilde{p}_{R,S}]
\square(3500,500)|alrm|/->`->`->`->/%
[\bullet`\bullet`\bullet`\bullet;p`1_{X\ox Y}`1_X`p]
\place(4250,1000)[=]
\square(4500,0)|alrb|/->`->`->`->/<500,1000>%
[\bullet`\bullet`\bullet`\bullet;p`R\ox S`R`p]
\morphism(4700,500)|a|<150,0>[`;\tilde{p}_{R,S}]
\efig$$

If all instances of $p$ in the sequence of diagrams
directly above are replaced by $r$ (with accompanying
changes of codomains), the equations
continue to hold and the two sequences of equations then jointly
establish, by Corollary \ref{!ness}, the first unitary 
coherence equation. Derivation of the second unitary coherence 
equation is similar.

Thus we have proved:
\thm\label{oxlax}
For a precartesian bicategory $\bB$, the data 
$(\ox,\widetilde\ox,\ox^\circ)$ constitute a lax 
functor $\bB\x\bB\ra\bB$.
\frp
\eth

\rmk\label{morelax}
Our main interest in the $\tilde p_{R,S}$ and the $\tilde r_{R,S}$
is in their role as components of the $p_{R,S}$ and the $r_{R,S}$,
which are the components of the pseudonatural transformations 
comprising the product projections $p\f P\la\ox\ra R\f r$, for the 
product pseudofunctor $\ox\f\bG\x\bG\ra\bG$. However we can say 
more. The naturality of the $\tilde p_{R,S}$ and the $\tilde r_{R,S}$
in $R$ and $S$ noted in Remark \ref{lxnat}, together with the equations 
{\em defining} $\ox^\circ$ and $\tilde\ox$, show that the 
$\tilde p_{R,S}$ and the $\tilde r_{R,S}$ provide lax naturality
2-cells making the $p_{X,Y}$ and the $r_{X,Y}$ components of
{\em lax natural} transformations $p\f P\la\ox\ra R\f r$, whose
domain is the lax functor 
$(\ox,\widetilde\ox,\ox^\circ)\f\bB\x\bB\ra\bB$.
It then follows that the $\tilde d_R$ provide lax naturality
2-cells making the $d_X$ components of a {\em lax natural}
transformation $d\f1_{\bB}\ra\ox D\f\bB\ra\bB$, 
where $D\f\bB\ra\bB\x\bB$ is the diagonal pseudofunctor.
\eth

\subsection{}\label{I}
Next, for a precartesian bicategory $\bB$, we describe a lax functor
$\one\ra\bB$, which amounts to giving an object of $\bB$ and
a monad on this. The object we take is $I$, the  terminal
object of $\mathrm{Map}(\bB)$.
The monad we take on the object $I$ has underlying arrow
$\top=\top_{I,I}$, the terminal object of $\bB(I,I)$.
The multiplication $\widetilde I\f\top\top\ra\top$ and the unit
$I^\circ\f1_I\ra\top$ are the unique 2-cells into the terminal 
object, which trivially satisfy the three monad equations. 

\prp\label{Ilax}
For a precartesian bicategory $\bB$,
the object $I$ of $\bB$, the arrow
$\top\f I\ra I$, and the 2-cells
$\widetilde I\f \top\top\ra\top$ and $I^\circ\f 1_I\ra\top$
constitute a lax functor $I\f\one\ra\bB$.
\frp
\eth

\rmk\label{Imorelax}
Further to Remark \ref{morelax} we note that the $\tilde t_R$
provide lax naturality 2-cells making the $t_X$ components of a 
{\em lax natural} transformation $t\f1_{\bB}\ra I!\f\bB\ra\bB$,
where $!\f\bB\ra\one$ is the unique such pseudofunctor.
\eth

\subsection{}\label{m}
We have the inclusion pseudofunctor $i\f\bM\ra\bB$. Given maps 
$f\f X\ra A$ and $g\f Y\ra B$ (in $\bM$ of course) we have, as
in \ref{prt}, (in $\bM$ and hence in $\bB$) the map 
$f\x g\f X\x Y\ra A\x B$ with invertible 2-cells
$p'_{f,g}\f p_{A,B}(f\x g)\ra f p_{X,Y}$ and
$r'_{f,g}\f r_{A,B}(f\x g)\ra g r_{X,Y}$. 
Since $f\ox g\f X\ox Y\ra A\ox B$ is the product of $f$ and $g$
in $\bG$, there is by Lemma \ref{eso} a unique 2-cell 
$m'_{f,g}\f f\x g\ra f\ox g$ for which the arrow
$(1,m'_{f,g},1)\f(X\x Y,f\x g,A\x B)\ra(X\ox Y,f\ox g,A\ox B)$
of $\bG$ satisfies the following two equations:
$$p_{f,g}(1,m'_{f,g},1)=(p_{X,Y},p'_{f,g},p_{A,B}),$$
$$r_{f,g}(1,m'_{f,g},1)=(r_{X,Y},r'_{f,g},r_{A,B}).$$

For the pseudofunctor $\x\f\bM\x\bM\ra\bM$ we can take
$1_X\x1_Y$ to be $1_{X\x Y}$ with $p'_{1,1}$ and $r'_{1,1}$ identities,
and then take $\x^{\circ}\f1_{X\x Y}\ra1_{X}\x1_{Y}$ 
to be again an identity;
while $\tilde{\x}\f(u\x v)(f\x g)\ra uf\x vg$ is the evident 2-cell.

For objects $X$ and $Y$ of $\bM$ we define 
$m_{X,Y}\f X\x Y\ra X\ox Y$ to be the identity.

\prp\label{laxnat}
For a precartesian bicategory $\bB$, 
the $m_{X,Y}$ and the $m'_{f,g}$ above
constitute a lax natural transformation
$m$ as in
$$\bfig
\square(0,0)/->`<-`<-`->/<1000,500>[\bB\x\bB`\bB`%
\bM\x\bM`\bM;\ox`i\x i`i`\x] 
\morphism(550,125)|m|<0,250>[`;m]
\efig$$
\eth
\prf
For each lax-naturality equation we appeal to Corollary \ref{!ness} 
For naturality in $f$ and $g$, we have 
$$p_{f,g}(1,m'_{f,g},1)(1,\phi\x\gamma,1)=%
(p_{X,Y},p'_{f,g},p_{A,B})(1,\phi\x\gamma,1)=%
(1,\phi,1)(p_{X,Y},p'_{f',g'},p_{A,B}),$$
where the first equality is by \ref{m} and the second is by
the pseudonaturality of $p$; while
$$p_{f,g}(1,\phi\ox\gamma,1)(1,m'_{f',g'},1)=%
(1,\phi,1)\tilde{p}_{f',g}(1,m'_{f',g'},1)=%
(1,\phi,1)(p_{X,Y},p'_{f',g'},p_{A,B}),$$
where the first equality is by \ref{aoxb} and the second is by the
equations defining the $m'_{f,g}$ in~\ref{m}.
Similarly, $(1,m'_{f,g},1)(1,\phi\x\gamma,1)$ and
$(1,\phi\ox\gamma,1)(1,m'_{f',g'},1)$ have the same composite with
$\tilde{r}_{f,g}$, so that Corollary \ref{!ness} gives 

$$\bfig
\morphism(0,0)|a|/{@{->}@/^4em/}/<750,0>[X\ox Y`A\ox B;f\ox g]
\morphism(0,0)|m|<750,0>[X\ox Y`A\ox B;f\x g]
\morphism(0,0)|b|/{@{->}@/_4em/}/<750,0>[X\ox Y`A\ox B;f'\x g']
\morphism(375,150)|l|<0,150>[`;m'_{f,g}]
\morphism(375,-250)|l|<0,150>[`;\phi\x\gamma]
\place(1000,400)[=]
\morphism(1250,0)|a|/{@{->}@/^4em/}/<750,0>[X\ox Y`A\ox B;f\ox g]
\morphism(1250,0)|m|<750,0>[X\ox Y`A\ox B;f'\ox g']
\morphism(1250,0)|b|/{@{->}@/_4em/}/<750,0>[X\ox Y`A\ox B;f'\x g']
\morphism(1625,150)|r|<0,150>[`;\phi\ox\gamma]
\morphism(1625,-250)|r|<0,150>[`;m'_{f',g'}]
\efig$$

For the nullary coherence condition, we have seen that
$\x^\circ$ and $p'_{1,1}$ can be taken to be identities; and then
$$p_{1,1}(1,m'_{1,1},1)(1,\x^\circ,1)=%
(p_{X,Y},p'_{1,1},p_{X,Y})=%
(p_{X,Y},1,p_{X,Y}),$$
where the first equality is by \ref{m}, while
$$p_{1,1}(1,\ox^\circ,1)=(p_{X,Y},1,p_{X,Y})$$ 
by \ref{lax}. Similarly with $r_{1,1}$ in place of $p_{1,1}$,
whereupon Corollary \ref{!ness} gives
$$\bfig
\morphism(0,0)|a|/{@{->}@/^4em/}/<750,0>[X\ox Y`X\ox Y;1_X\ox1_Y]
\morphism(0,0)|m|<750,0>[X\ox Y`X\ox Y;1_X\x1_Y]
\morphism(0,0)|b|/{@{->}@/_4em/}/<750,0>[X\ox Y`X\ox Y;1_{X\ox Y}]
\morphism(375,100)|m|<0,250>[`;m'_{1_X,1_Y}]
\morphism(375,-300)|m|<0,200>[`;\x^{\circ}]
\place(1000,400)[=]
\morphism(1250,0)|a|/{@{->}@/^4em/}/<750,0>[X\ox Y`X\ox Y;1_X\ox1_Y]
\morphism(1250,0)|b|/{@{->}@/_4em/}/<750,0>[X\ox Y`X\ox Y;1_{X\ox Y}]
\morphism(1625,-100)|m|<0,250>[`;\ox^{\circ}]
\efig$$

Finally, for the binary coherence equation we require, 
for $f\f X\ra A$, $g\f Y\ra B$, $u\f A\ra L$, and $v\f B\ra M$ 
in $\bM$,
$$\bfig
\morphism(0,0)|a|/{@{->}@/^2em/}/<500,0>[X\ox Y`A\ox B;f\ox g]
\morphism(0,0)|b|/{@{->}@/_2em/}/<500,0>[X\ox Y`A\ox B;f\x g]
\morphism(250,-100)|m|<0,250>[`;m'_{f,g}]
\morphism(500,0)|a|/{@{->}@/^2em/}/<500,0>[A\ox B`L\ox M;u\ox v]
\morphism(500,0)|b|/{@{->}@/_2em/}/<500,0>[A\ox B`L\ox M;u\x v]
\morphism(750,-100)|m|<0,250>[`;m'_{u,v}]
\morphism(0,0)|a|/{@{->}@/^5em/}/<1050,0>%
[X\ox Y`\phantom{A\ox B};uf\ox vg]
\morphism(500,200)|m|<0,225>[`;\widetilde\ox]
\place(-250,400)[=]
\morphism(-1500,0)|a|/{@{->}@/^4em/}/<1000,0>[X\ox Y`L\ox M;uf\ox fg]
\morphism(-1500,0)|m|<1000,0>[X\ox Y`L\ox M;uf\x vg]
\morphism(-1500,0)|b|/{@{->}@/_1em/}/<500,-500>[X\ox Y`A\ox B;f\x g]
\morphism(-1000,-500)|b|/{@{->}@/_1em/}/<500,500>[A\ox B`L\ox M;u\x v]
\morphism(-1000,-350)|m|<0,200>[`;\widetilde\x]
\morphism(-1000,100)|m|<0,250>[`;m'_{uf,vg}]
\efig$$
Now the composite of the right side with $p_{uf,vg}$ equals
by \ref{lax} the vertical pasting of $p_{f,g}(1,m'_{f,g},1)$
and $p_{u,v}(1,m'_{u,v},1)$ along $p_{A,B}$, and thus by \ref{m}
equals the vertical pasting along $p_{A,B}$ of 
$(p_{X,Y},p'_{f,g},p_{A,B})$
and $(p_{A,B},p'_{u,v},p_{L,M})$---which, 
because $p$ is pseudonatural,
equals the composite of $(p_{X,Y},p'_{uf,vg},p_{L,M})$ 
with $\tilde{\x}$.
But by \ref{m} this is also equal to the composite of the left side 
with $p_{uf,vg}$. In the same way the two sides have equal composites 
with $r_{uf,vg}$, whence they are equal by Corollary \ref{!ness}. 
\frp

\subsection{}\label{2.19}
Consider objects $R\f X\ra A$ and $S\f Y\ra B$ of $\bG$, 
along with their
product $R\ox S$ and its projections (in their primary forms)
$(p_{X,Y},\tilde{p}_{R,S},p_{A,B})$ and 
$(r_{X,Y},\tilde{r}_{R,S},r_{A,B})$. Consider also
maps $f\f L\ra X$ and $g\f M\ra Y$, along with the invertible 2-cells
$p'_{f,g}:p_{X,Y}(f\x g)\ra fp_{L,M}$ and 
$r'_{f,g}:r_{X,Y}(f\x g)\ra gr_{L,M}$ of \ref{prt}. 
Now form the vertical pasting of $(p_{X,Y},\tilde{p}_{R,S},p_{A,B})$ 
and $(p_{L,M},p'_{f,g},p_{X,Y})$ along $p_{X,Y}$ and similarly the 
vertical pasting of $(r_{X,Y},\tilde{r}_{R,S},r_{A,B})$ and 
$(r_{L,M},r'_{f,g},r_{X,Y})$ along $r_{X,Y}$.
These pastings constitute arrows in $\bG$ from $(R\ox S)(f\x g)$ to 
$Rf$ and to $Sg$, and thus determine an arrow in $\bG$ from 
$(R\ox S)(f\x g)$ to $Rf\ox Sg$. 

In fact this arrow is an isomorphism. To see this, observe that the
secondary forms of these pastings are made by pasting 
$(p_{L,M},p'_{f,g},p_{X,Y})$ and $(r_{L,M},r'_{f,g},r_{X,Y})$ to the
secondary forms of $(p_{X,Y},\tilde{p}_{R,S},p_{A,B})$ and 
$(r_{X,Y},\tilde{r}_{R,S},r_{A,B})$, which are 
$(p_{X,Y},\pi,p_{A,B})$ and
$(r_{X,Y},\rho,r_{A,B})$ where $\pi$ and $\rho$ are the projections
from $R\ox S = p^*Rp\wedge r^*Sr$. Since {\em precomposition} with the
map $f\x g$, being a right adjoint, preserves products, the projections
$\pi(f\x g)$ and $\rho(f\x g)$ express $(R\ox S)(f\x g)$ as a product
$p^*Rp(f\x g)\wedge r^*Sr(f\x g)$. Then, because the 2-cells
$p'_{f,g}\f p_{X,Y}(f\x g)\ra fp_{L,M}$ and 
$r'_{f,g}\f r_{X,Y}(f\x g)\ra gr_{L,M}$ 
are invertible, $(R\ox S)(f\x g)$ is also a product 
$p^*Rfp\wedge r^*Sgr$ with projections the vertical pastings of
$(p_{X,Y},\pi,p_{A,B})$ to $(p_{L,M},p'_{f,g},p_{X,Y})$ and of
$(r_{X,Y},\rho,r_{A,B})$ to $(r_{L,M},r'_{f,g},r_{X,Y})$. This 
completes the proof. Similarly, for maps $u\f Z\ra A$ and $v\f W\ra B$,
we can construct an isomorphism 
$(u\x v)^*(R\ox S){\s \to^{\simeq}}\,u^*R\ox v^*S$. (These isomorphisms
should be seen in the light of the ${\bf L}$-homomorphisms
of~(2.15) of \cite{ckvw}.)

If we apply the first result above with $A=X$, $B=Y$, $R=1_X$ and
$S=1_Y$, we have an isomorphism $\xi\f (1_X\ox 1_Y)(f\x g)\ra f\ox g$,
whose projections are the pasting of $(p_{L,M},p'_{f,g},p_{X,Y})$
to $(p_{X,Y},p_{1,1},p_{X,Y})$ and the same with $r$ in place of $p$.
Composing $\xi$ with $\ox^\circ(f\x g)$ gives by 2.13 a 2-cell
$f\x g\ra f\ox g$ whose projections are the pasting of
$(p_{L,M},p'_{f,g},p_{X,Y})$ to an identity, and the same with $r$ in
place of $p$. By \ref{m}, therefore, this arrow 
$\xi.\ox^\circ(f\x g)\f f\x g\ra f\ox g$ is $m'_{f,g}$.

\prp\label{laxnatI}
For a precartesian bicategory $\bB$, there is a lax natural 
transformation $u$ as in
$$\bfig
\Dtriangle (1000,0)|arb|/<-`<- `->/<500,250>%
[\bB`\one`\bM;i`I`1]
\morphism(1150,200)|m|<0,150>[`;u]
\efig$$
whose component $u$ at the unique object of $\one$ is
$1_I\f I\ra I$ and whose lax naturality 2-cell at the
unique arrow of $\one$ is $\tau\f1\ra\top$, the
lax naturality here being equivalent to the two equations: 
$$\bfig
\morphism(0,0)|a|/{@{->}@/^2em/}/<500,0>[I`I;\top]
\morphism(0,0)|b|/{@{->}@/_2em/}/<500,0>[I`I;1_I]
\morphism(250,-100)|m|<0,225>[`;\tau]
\place(750,0)[=]
\morphism(1000,0)|a|/{@{->}@/^2em/}/<500,0>[I`I;\top]
\morphism(1000,0)|b|/{@{->}@/_2em/}/<500,0>[I`I;1_I]
\morphism(1250,-100)|m|<0,225>[`;I^\circ]
\efig$$
and
$$\bfig
\morphism(0,0)|a|/{@{->}@/^3em/}/<500,0>[I`I;\top]
\morphism(0,0)|m|<500,0>[I`I;1_I]
\morphism(0,0)|b|/{@{->}@/_3em/}/<500,0>[I`I;1_I1_I]
\morphism(250,50)|m|<0,200>[`;\tau]
\place(750,0)[=]
\morphism(1000,0)|a|/{@{->}@/^2em/}/<500,0>[I`I;\top]
\morphism(1000,0)|b|/{@{->}@/_2em/}/<500,0>[I`I;1_I]
\morphism(1250,-100)|m|<0,225>[`;\tau]
\morphism(1500,0)|a|/{@{->}@/^2em/}/<500,0>[I`I;\top]
\morphism(1500,0)|b|/{@{->}@/_2em/}/<500,0>[I`I;1_I]
\morphism(1750,-100)|m|<0,225>[`;\tau]
\morphism(1000,0)|a|/{@{->}@/^5em/}/<1050,0>[I`\phantom{I};\top]
\morphism(1500,200)|m|<0,225>[`;\widetilde I]
\efig$$
which hold automatically since $\top$ is the terminal
object of $\bB(I,I)$.
\frp
\eth

\prp\label{psnat}
For a precartesian bicategory $\bB$, the lax natural transformation 
$$m:i.\x\ra\ox.(i\x i):\bM\x\bM\ra\bB$$
is a pseudonatural transformation if and only if
the constraints $\ox^{\circ}$ are invertible, and then
$m$ itself is invertible. Again, the lax transformation
$$u : i.1\ra I : \one\ra\bB$$
is a pseudonatural transformation if and only if
the constraint $I^\circ\f 1_I\ra\top$ is invertible, and then
$\tilde{I}$ and $u$ itself are invertible. 
\eth
\prf
The lax natural $m$ is pseudonatural when the 2-cells $m'_{f,g}$
are invertible.
By the last paragraph of \ref{2.19}, each $m'_{f,g}$ is invertible
when $\ox^{\circ}$ is invertible; on the other hand 
$\ox^{\circ}$ is invertible if $m'_{1_X,1_Y}$ is invertible
by the second equation of Proposition \ref{laxnat}. Then $m$ is invertible, because
each $m_{X,Y}$ is an identity. The truth of the assertions
about $u$ is immediate from Proposition \ref{laxnatI}.
\frp

\section{Cartesian Bicategories}\label{cartbicat}
\dfn\label{main} A precartesian bicategory $\bB$ is said to be
{\em cartesian\/} when the
$$\bB\x\bB\to^{\ox}\bB\toleft^{I}\one$$
are pseudofunctors, meaning that $\widetilde\ox$, $\ox^\circ$,
$I^\circ$ (and hence $\widetilde I$) are invertible.
\eth

Since pseudofunctors carry adjunctions to adjunctions, it should
be noted that in a cartesian bicategory, the map $f\x g$ arising
from maps $f$ and $g$ has adjunction data
$$\eta_f\ox\eta_g,\epsilon_f\ox\epsilon_g\f f\x g\laj f^*\ox g^*$$ 

\prp\label{unifier} For a cartesian bicategory $\bB$,
the pseudofunctors
$$\bB\x\bB\to^{\ox}\bB\toleft^{I}\one$$
\begin{enumerate}
\item[1)] restrict to $\bM$ giving a right
adjoint to 
$$\bM\x\bM\toleft^{\Delta}\bM$$
with unit $d$ and a right adjoint to $$\bM\to^{!}\one$$
with unit $t$;
\item[2)] the composites 
$$\bB(X,Y)\x\bB(X,Y)\to^{\ox}\bB(X\ox X,Y\ox Y)\to^{\bB(d_X,d_Y^*)}%
\bB(X,Y)$$
and
$$\bB(X,Y)\toleft^{\bB(t_X,t_Y^*)}\bB(I,I)\toleft^{\top}\one$$
provide right adjoints to
$$\bB(X,Y)\x\bB(X,Y)\toleft^{\Delta}\bB(X,Y)\quad\mbox{and}\quad%
\bB(X,Y)\to^{!}\one$$
\end{enumerate}
Moreover, for any other pair of pseudofunctors 
$\bB\x\bB\to^{\ox'}\bB\toleft^{I'}\one$, if they satisfy~1) 
and~2) then $\ox'\cong \ox$ and $I'\cong I$ (as pseudofunctors). 
Thus, a cartesian bicategory
can be described alternatively as a bicategory $\bB$ together
with pseudofunctors $\ox$ and $I$, as under consideration, 
which satisfy 1) and 2).
\eth
\prf
We have seen that $\ox$ and $I$ satisfy the conditions of the 
proposition. Only their essential uniqueness remains to be shown.
If $\ox'$ and $I'$ satisfy {\em 1)\/} 
then their satisfaction of {\em 2)\/} 
is equivalent to their providing right adjoints to 
$$\bG\x\bG\toleft^{\Delta}\bG\to^!\one$$
commuting via $\partial_0$ and $\partial_1$ with the
corresponding adjoints for $\bM$. By essential uniqueness
of products we have diagrams
$$\bfig
\square(0,0)[X\ox' Y`X\ox Y`A\ox' B`A\ox B;%
k_{X,Y}`R\ox' S`R\ox S`k_{A,B}]
\morphism(125,250)|a|<250,0>[`;\tilde k_{R,S}]
\square(1250,0)[I'`I`I'`I;l_*`\top`\top'`l_*]
\morphism(1375,255)|a|<250,0>[`;\tilde l_{1_*}]
\efig$$
which are equivalences in $\bG$. By Proposition \ref{Gequiv}
the $k_{X,Y}$ and $l_*$ are equivalences in $\bM$ and
the $\tilde k_{R,S}$ and $\tilde l_{1_*}$ are isomorphisms. 
These provide the components for pseudonatural
equivalences $k\f\ox'\ra\ox$ and $l\f I'\ra I$.
\frp

\rmk Since a terminal object in a (mere) category is unique 
to within unique isomorphism and we have $I^\circ\f1_I\ra\top$
invertible in a cartesian bicategory we may as well as assume,
in a cartesian bicategory, that we have chosen $\top=1_I$.
\eth

\subsection{}\label{herewego} 
We turn now to further analysis of the pseudofunctors
$$\bB\x\bB\to^{\ox}\bB\toleft^{I}\one$$
Since $\ox$ is binary product in $\bM$ we have the pseudonatural
(adjoint) equivalences
$$a_{X,Y,Z}\f(X\ox Y)\ox Z\to X\ox(Y\ox Z)\qquad%
l\f I\ox X\to X\qquad r\f X\ox I\to X$$
$$\mbox{and}\qquad s\f X\ox Y\to Y\ox X$$
in $\bM$ as constructed and studied in Section \ref{finprod}. 
However, we can say more:
  
\prp\label{psnatcon}
The equivalence maps $a\f(X\ox Y)\ox Z\ra X\ox (Y\ox Z)$,
$l\f I\ox X\ra X$, $r\f X\ox I\ra X$, and $s\f X\ox Y\ra Y\ox X$
extend to pseudonatural equivalences between the relevant
$\bB$-valued functors.
\eth
\prf
Because $\bG$ has finite products which commute via $\partial_0$
and $\partial_1$ with those of $\bM$ we have equivalences
$$\bfig
\scalefactor{0.94}
\square(0,0)<800,500>%
[(X\ox Y)\ox Z`X\ox(Y\ox Z)`(A\ox B)\ox C`A\ox(B\ox C);%
a`(R\ox S)\ox T`R\ox(S\ox T)`a]
\morphism(225,250)|m|<350,0>[`;\tilde a_{R,S,T}]
\square(1400,0)[I\ox X`X`I\ox A`A;l`I\ox R`R`l]
\morphism(1525,250)|m|<250,0>[`;\tilde l_{R}]
\square(2200,0)[X\ox I`X`A\ox I`A;r`R\ox I`R`r]
\morphism(2325,250)|m|<250,0>[`;\tilde r_{R}]
\square(3000,0)[X\ox Y`Y\ox X`A\ox B`B\ox A;%
s`R\ox S`S\ox R`s]
\morphism(3125,250)|m|<250,0>[`;\tilde s_{R,S}]
\efig$$
in $\bG$. {}From Proposition \ref{Gequiv} we have that the 
$\tilde a_{R,S,T}$, $\tilde l_{R}$, $\tilde r_{R}$, and 
$\tilde s_{R,S}$ are invertible so that the squares
above can be seen as providing the data for 
pseudonatural transformations
between $\bB$-valued functors. We omit verification
of the pseudonaturality equations.
\frp
   
\thm\label{smbicat}
The data $(\bB,\ox,I,a,l,r,s)$ extend to provide a
symmetric monoidal bicategory.
\eth
\prf
By Theorem \ref{finprodimpmonl} the bicategories $\bM$ and $\bG$
are symmetric monoidal bicategories. To construct the modifications
$\pi$, $\mu$, $\lambda$, and $\rho$ of \cite{gps} observe that we
have invertible 2-cells $\pi_{X,Y,Z,W}$ as in the first pentagonal 
region below (in which we have abbreviated $\ox$ by juxtaposition) 
because $\bM$ is monoidal. Since $\bM$ is a subbicategory
of $\bB$ the $\pi$ are also 2-cells in $\bB$. For these to constitute
a $\bB$-valued modification we require the following equality of
pastings.
$$\bfig
\scalefactor{0.75}
\morphism(0,0)|l|<500,-500>[((XY)Z)W`(X(YZ))W;aW]
\morphism(500,-500)|b|<1000,0>[(X(YZ))W`X((YZ)W));a]
\morphism(1500,-500)|r|<500,500>[X((YZ)W))`X(Y(ZW));Xa]
\morphism(0,0)<1000,500>[((XY)Z)W`(XY)(ZW);a]
\morphism(1000,500)<1000,-500>[(XY)(ZW)`X(Y(ZW));a]
\morphism(1000,-175)|m|<0,250>[`;\pi_{X,Y,Z,W}]
\morphism(0,-1000)|l|<500,-500>[((AB)C)D`(A(BC))D;aD]
\morphism(500,-1500)|b|<1000,0>[(A(BC))D`A((BC)D));a]
\morphism(1500,-1500)|r|<500,500>[A((BC)D))`A(B(CD));Aa]
\morphism(0,0)|l|<0,-1000>[((XY)Z)W`((AB)C)D;((RS)T)U]
\morphism(500,-500)<0,-1000>[(X(YZ))W`(A(BC))D;]
\morphism(1500,-500)<0,-1000>[X((YZ)W))`A((BC)D));]
\morphism(2000,0)<0,-1000>[X(Y(ZW))`A(B(CD));]
\morphism(250,-900)|m|<0,250>[`;\tilde\alpha U]
\morphism(1000,-1150)|m|<0,250>[`;\tilde\alpha]
\morphism(1750,-900)|m|<0,250>[`;R\tilde\alpha]
\place(1000,-1750)[=]
\morphism(0,-2500)<1000,500>[((XY)Z)W`(XY)(ZW);a]
\morphism(1000,-2000)<1000,-500>[(XY)(ZW)`X(Y(ZW));a]
\morphism(0,-3500)|l|<500,-500>[((AB)C)D`(A(BC))D;aD]
\morphism(500,-4000)|b|<1000,0>[(A(BC))D`A((BC)D));a]
\morphism(1500,-4000)|r|<500,500>[A((BC)D))`A(B(CD));Aa]
\morphism(0,-3500)<1000,500>[((AB)C)D`(AB)(CD);a]
\morphism(1000,-3000)<1000,-500>[(AB)(CD)`A(B(CD));a]
\morphism(0,-2500)<0,-1000>[((XY)Z)W`((AB)C)D;((RS)T)U]
\morphism(1000,-2000)<0,-1000>[(XY)(ZW)`(AB)(CD);]
\morphism(2000,-2500)<0,-1000>[X(Y(ZW))`A(B(CD));]
\morphism(1000,-3675)|m|<0,250>[`;\pi_{A,B,C,D}]
\morphism(500,-2850)|m|<0,250>[`;\tilde\alpha]
\morphism(1500,-2850)|m|<0,250>[`;\tilde\alpha]
\efig$$
However, this equation follows simply because the
pair $(\pi_{W,X,Y,Z},\pi_{A,B,C,D})$ is an invertible
2-cell in the monoidal bicategory $\bG$. Finally, to
see for example that the $\pi$ satisfy the non-abelian
cocycle condition in \cite{gps} observe that the
pentagons and squares in that diagram have boundaries
composed entirely of maps so that the condition is
satisfied because $\bM$ is a monoidal bicategory. The other
symmetric monoidal bicategory data is dealt with in a 
similar manner.
\frp

As we stressed earlier the 2-cells in the product 
projections displayed in
(\ref{projexionsa}) are, in general, by no means invertible
in $\bB$. However, for those of the special form $\tilde p_{R,1_Y}$ 
or $\tilde r_{1_X,S}$ we have:

\prp\label{spiso} In a cartesian bicategory, for any arrow
$R\f X\ra A$ and any object $Y$, the 2-cells
$$\bfig
\square(0,0)[X\ox Y`X`A\ox Y`A;p_{X,Y}`R\ox 1_Y`R`p_{A,Y}]
\morphism(125,250)|m|<300,0>[`;\tilde p_{R,1_Y}]
\place(1000,500)[\mbox{and}]
\square(1500,0)[X\ox Y`X\ox I`A\ox Y`A\ox I;%
1_X\ox t_Y`R\ox 1_Y`R\ox \top`1_A\ox t_Y]
\morphism(1600,250)|m|<350,0>[`;1_R\ox\tilde t_{1_Y}]
\efig$$
are invertible. Similarly, for any object $X$ and any arrow
$S\f Y\ra B$, the 2-cell $r_{1_X,S}$ is invertible.
\eth
\prf In any bicategory with finite products,
the equivalence $X\ra X\x I$ identifies 
$p_{X,Y}\f X\x Y\ra X$ and $1_X\x t_Y\f X\x Y\ra X\x I$
to within isomorphism. 
In particular, this applies to $\bG$ where the inverse
of the relevant equivalence is
$$\bfig
\square(0,0)[X\ox I`X`A\ox I`A;p_{X,I}`R\ox 1_I`R`p_{A,I}]
\morphism(125,250)|m|<300,0>[`;\tilde p_{R,1_I}]
\efig$$
in which, by Proposition \ref{Gequiv}, the 2-cell $\tilde p_{R,1_I}$
is invertible. When the equivalence square above is pasted to the
square on the right in the statement the result is the square on the
left in the statement. In the square on the right of the statement
observe first that since $\top =1_I$ we have $\tilde t_{1_Y}=1_{t_Y}$
by uniqueness of such arrows to $\top\f I\ra I$ in $\bG$.
Now $1_R\ox\tilde t_{1_Y}=1_R\ox 1_{t_Y}$ is invertible being a 
$\ox$-product of invertible 2-cells. (We note for clarity though that 
our notation here has $1_R\ox 1_{t_Y}$ being the $\ox$-product
of the identity $1_AR \ra R1_X$ and the identity $t_Y1_Y\ra 1_It_Y$
by the assumed normality of $\bB$. As a paste composite of 
invertibles, $\tilde p_{R,1_Y}$ is also invertible. Similarly,
since $1_{t_X}\ox1_S$ is invertible, $\tilde r_{1_X,S}$ is invertible.
\frp

A similar result holds for the mates of the $\tilde p_{R,1_Y}$ and
the $\tilde r_{1_X,S}$.

\prp\label{preBeck} In a cartesian bicategory, for any arrow
$R\f X\ra A$ and any object $Y$, the mate of $\tilde p_{R,1_Y}$ as in
$$\bfig
\square(1000,0)/<-`->`->`<-/%
[X\ox Y`X`A\ox Y`A;p^*_{X,Y}`R\ox1_Y`R`p^*_{A,Y}]
\morphism(1125,250)|m|<300,0>[`;\tilde p^*_{R,1_Y}]
\efig$$
is invertible. Similarly, for any arrow $S\f Y\ra B$ and any object $X$,
the mate $\tilde r^*_{1_X,S}$ of $\tilde r_{1_X,S}$ is invertible.
\eth
\prf
As in Proposition \ref{spiso} the task reduces to showing that
the mate of $1_R\ox\tilde t_{1_Y}= 1_R\ox1_{t_Y}$ with respect to the adjunctions
$1_X\ox t_Y\laj(1_X\ox t_Y)^*$ and $1_A\ox t_Y\laj(1_A\ox t_Y)^*$
is invertible. Since $\ox$ is a psedofunctor it preserves adjunctions
and the mate operation. The mate at issue is easily seen to be
$1_R\ox1_{t^*_Y}$ which being a $\ox$-product of invertibles is
invertible.
\frp

\prp\label{strange}
For $X$ and $Y$ in a cartesian bicategory there is a natural
isomorphism
$$\bfig
\Atriangle[\bB(X,I)\x\bB(I,Y)`\bB(X\ox I,I\ox Y)`\bB(X,Y);%
\ox`\circ`\bB(p^*,r)]
\place(500,150)[\iso]
\efig$$
where we have used $\circ$ to denote composition in $\bB$. Moreover,
the 2-cells
$$\bfig
\square(500,0)|almm|/<-`->`->`<-/[X`X`I`I;1`R`R`1]
\square(500,-500)|mlmb|/<-`->`->`<-/[I`I`I`Y;1`\top`S`t]
\morphism(875,-250)|m|<-250,0>[`;\tau]
\morphism(500,500)|l|/{@{->}@/_2em/}/<0,-1000>[X`I;R]
\square(1000,0)|amrm|[X`I`I`I;t`R`\top`1]
\morphism(1125,250)|m|<250,0>[`;\tau]
\square(1000,-500)|mmrb|[I`I`Y`Y;1`S`S`1]
\morphism(1500,500)|r|/{@{->}@/^2em/}/<0,-1000>[I`Y;S]
\efig$$
provide product projections for $SR\f X\ra Y$ seen
as a product of $R\f X\ra I$ and $S\f I\ra Y$ in $\bG$.
\eth
\prf
We have
$$\bfig
\morphism(0,1000)|l|<0,-500>[X`X\ox I;p^*]
\morphism(0,1000)|a|<500,-500>[X`I;R]
\morphism(500,500)|l|<0,-500>[I`I\ox I;p^*]
\Dtriangle(0,-500)|lll|[X\ox I`I\ox I`I\ox Y;R\ox S`R\ox I`I\ox S]
\morphism(0,-500)|l|<0,-500>[I\ox Y`Y;r]
\morphism(0,-500)|l|<0,-500>[I\ox Y`Y;r]
\morphism(500,0)|l|<0,-500>[I\ox I`I;r]
\morphism(500,-500)|r|<-500,-500>[I`Y;S]
\morphism(500,500)|r|/{@{->}@/^3em/}/<0,-1000>[I`I;\top]
\morphism(300,0)|m|<-150,0>[`;\tilde\ox]
\morphism(250,350)|a|<150,0>[`;\tilde p^*_{R,1}]
\morphism(250,-350)|b|<150,0>[`;\tilde r_{1,S}]
\morphism(625,0)|m|<150,0>[`;\tau]
\efig$$
in which all the 2-cells are invertible, those in the
parallelograms being so by Propositions~\ref{preBeck}
and \ref{spiso}. {}From this explicit description the second
clause of the statement follows easily.
\frp

\end{document}